%% file: main.tex
\documentclass[11pt,a4paper,USenglish]{article}
\usepackage{template}
\usepackage{fix-cm}
\usepackage{document_style}
\usepackage{indentfirst}

\title{A Battle-Lemari\'e Frame Characterization \\ of Besov and Triebel-Lizorkin Spaces}
\author{Andrew Haar\inst{1}}
\address{andrew.haar674@student.cuni.cz}

\begin{document}

\maketitle

{\centerline{\emph{This article is dedicated to Prof. Hans Triebel on the Occasion of his 90th Birthday.}}}

\begin{abstract}
In this paper, we investigate a spline frame generated by oversampling against the well-known Battle-Lemari\'e wavelet system of nonnegative integer order, $n$. We establish a characterization of the Besov and Triebel-Lizorkin (quasi-) norms for the smoothness parameter up to $s < n+1$, which includes values of $s$ where the Battle-Lemari\'e system no longer provides an unconditional basis; we, additionally, prove a result for the endpoint case $s=n+1$. This builds off of earlier work by G. Garrig\'os, A. Seeger, and T. Ullrich, where they proved the case $n=0$, i.e. that of the Haar wavelet, and work of R. Srivastava, where she gave a necessary range for the Battle-Lemari\'e system to give an unconditional basis of the Triebel-Lizorkin spaces.\blankfootnote{2020 Mathematics Subject Classification: 46E35, 46B15, 42C40, 42B35.}\blankfootnote{\emph{Keywords}: Battle-Lemari\'e system, wavelet, spline, frame, Besov spaces, Triebel-Lizorkin spaces.}
\end{abstract}

\import{Sections}{introduction.tex}

\import{Sections}{preliminaries.tex}

\import{Sections}{mainResultProof.tex}

\import{Sections}{edgeCase.tex}

\printbibliography

\end{document}

%% file: Sections/introduction.tex
\section{Introduction}

The Battle-Lemari\'e wavelet system has been an important tool for many years because of its remarkable ability to characterize function spaces, particularly those measuring smoothness \cite{srivastava2023orthogonal,triebel2010bases,ushakova2018localisation}, even those that also have weights \cite{ushakova2021spline}. These types of function space decompositions through orthonormal wavelet systems have found their place in the literature not only due to their orthonormality, but also because the wavelets have vanishing moments (property (IV) in Section \ref{subsec:wavelets}); for these reasons, this system is particularly well-suited to study asymptotic decay rates for multiresolution approximations \cite{mallat2002theory,mallat1989multiresolution}, denoising from digital images \cite{mudugamuwa2012battle}, entropy numbers \cite{nasyrova2018wavelet}, etc. 

We will expand upon some known characterizations of the Besov and Triebel-Lizorkin spaces by creating a frame defined by oversampling against the Battle-Lemari\'e system; the novelty here is that we can widen the range of smoothness that we can measure past the point at which the Battle-Lemari\'e system ceases to give an unconditional basis of these function spaces. Our work generalizes a result from \cite{garrigos2023haar}, as suggested in \cite[Remark 1.7.3]{garrigos2023haar}, in which the authors demonstrated the three theorems that we will show, but for the order zero Battle-Lemari\'e system, i.e. the well-known Haar wavelet system \eqref{eq:HaarSystem}.

Our object of concern throughout this paper will be the so-called Battle-Lemari\'e orthogonal spline system, investigated by both G. Battle \cite{battle1987block} and P. Lemari\'e \cite{lemarie1988ondelettes} in the late 1980s. The Battle-Lemari\'e scaling function of order $n\in\nN_0 = \nN\cup \{0\}$, call it $\Psi_n$, and the corresponding wavelet, $\psi_n$, are real-valued, square integrable functions that give rise to an orthonormal basis of $L_2(\nR)$, the space of square integrable functions, in the following way. Write
\begin{align*}
    \psi_{n;j,\mu}(x) = \psi_n(2^j x - \mu) , \quad \psi_{n;-1,\mu}(x) = \sqrt{2} \Psi_n(x-\mu)
\end{align*}
for $j\in\nN_0$ and $\mu\in\nZ$. Then
\begin{align*}
    \cW_n = \{ 2^{j/2} \psi_{n;j,\mu} : j\in\nN_0\cup\{-1\} , \mu\in\nZ \}
\end{align*}
constitutes an orthonormal basis of $L_2(\nR)$. We will discuss further properties of these spline functions in Section \ref{subsec:wavelets}.

At this point, we fix $n\in\nN_0$ and cease to clutter the subscript of our wavelet and scaling function with it, i.e. we write $\psi$ and $\Psi$ rather than $\psi_n$ and $\Psi_n$. There are actually many Battle-Lemari\'e wavelet systems (see Definition \ref{def:bl}); at any given moment, $\psi,\Psi$ will be fixed.

We remark that, except when $n=0$, where we get the Haar system, these splines do not have compact support, but they do concentrate their masses around $0$, which means that their scaled versions, $\psi_{j,\mu}$, concentrate their masses near the following dyadic intervals: for $j\geqslant 0$ and $\mu\in\nZ$,
\begin{align}
    I_{j,\mu} = [2^{-j}\mu,2^{-j}(\mu+1)) . \label{eq:dyadicInterval}
\end{align}
We define $I_{-1,\mu} = [\mu,\mu+1)$ for $\mu\in\nZ$ since $\psi_{-1,\mu}(x) = \sqrt{2} \Psi(x-\mu)$ and \emph{not} $\sqrt{2}\Psi(2^{-1}x - \mu)$.

It is natural to ask from this standpoint if these spline wavelets can be used to characterize function spaces other than $L_2(\nR)$; the answer is a resounding yes! We will focus on the classic Besov and Triebel-Lizorkin spaces, denoted $B^s_{p,q}(\nR)$ and $F^s_{p,q}(\nR)$, respectively; see Section \ref{subsec:spaces} for definitions. 

Suppose, for now, that $1 < p,q < \infty$; then with our spline system, $\cW_n$, we get an unconditional basis of $B^s_{p,q}(\nR)$ for
\begin{align}
    -\frac{1}{p'} - n < s < n + \frac{1}{p} , \label{eq:sBesovIntro}
\end{align}
and for $F^s_{p,q}(\nR)$ with
\begin{align}
    \max\left\{ - \frac{1}{p'} , - \frac{1}{q'} \right\} - n < s < n + \min\left\{ \frac{1}{p} , \frac{1}{q} \right\} . \label{eq:sTriebelIntro}
\end{align}
By $p',q'$, we mean H\"older conjugates (see the Notation section). The results above are well-known and can be found in \cite[Theorems 2.46 and 2.49]{triebel2010bases}, see also Theorem \ref{thm:SplineCharac} below. As was proven in \cite{seeger2017haar} for the case $n=0$ and in \cite{srivastava2023orthogonal} for the general case, the given ranges of $s$ are also necessary for $\cW_n$ to give an unconditional basis in these spaces when $1<p,q<\infty$.

We are now prepared to define our frame. Recall that, for $H$, a Hilbert space, we say $(g_j)_j\subset H$ is a frame if there exist constants, $A,B>0$, such that 
\begin{align*}
    A \|x\|_H^2 \leqslant \sum_j |(x,g_j)_H|^2 \leqslant B \|x\|_H^2
\end{align*}
for every $x\in H$. Of course, $(\cdot,\cdot)_H$ and $\|\cdot\|_H$ refer to the inner product and norm on $H$, respectively.

Write
\begin{align}
    \widetilde{\psi}_{j,\mu}(x) = \psi\left(2^j x - \frac{\mu}{2}\right) \label{eq:shiftedSpline}
\end{align}
for $j\in\nN_0$ and $\mu\in\nZ$. For the base scale,
\begin{align}
    \widetilde{\psi}_{-1,\mu}(x) = \Psi(x-\mu) , \quad \mu\in\nZ . \label{eq:shiftedSpline-1}
\end{align}
Define
\begin{align*}
    \cW_n^{\mathrm{ext}} = \{\widetilde{\psi}_{-1,\mu} : \mu\in\nZ\} \cup \{2^{j/2}\widetilde{\psi}_{j,\mu} : j\in\nN_0, \mu\in\nZ\} ;
\end{align*}
this gives a frame for $L_2(\nR)$. Indeed, if we were to include $\{\Psi(x-\frac{2\mu + 1}{2}) : \mu\in\nZ\}$ in $\cW_n^{\mathrm{ext}}$ then it would be the union of two orthonormal bases, making it a frame. Obviously, $\cW_n^{\mathrm{ext}}$ is still a frame despite the omission of this extra set.

The corresponding \emph{shifted spline coefficients} will be
\begin{align}
    \sfk_{j,\mu}(f) = 2^j |(f,\widetilde{\psi}_{j,2\mu})| + 2^j |(f,\widetilde{\psi}_{j,2\mu+1})| \label{eq:shiftedSplineCoef}
\end{align}
for $j\in\nN_0$ and $\mu\in\nZ$ along with
\begin{align}
    \sfk_{-1,\mu}(f) = |(f,\widetilde{\psi}_{-1,\mu})| \label{eq:shiftedSplinej=-1}
\end{align}
also for $\mu\in\nZ$. 

The notation $(\cdot,\cdot)$ should be intuitively thought of as a pairing of a distribution with a function. However, because $\psi$ is not infinitely differentiable (property (I) in Section \ref{subsec:wavelets}), the pairing is not defined \emph{a priori} for all $f\in\cS'(\nR)$. We will give a definition of this pairing in Section \ref{sec:duality} that is defined for all $f\in\cS'(\nR)$, although we will be using $\cB\subset \cS'$, defined by
\begin{align*}
    \cB = B^{-(n+1)}_{\infty,1}(\nR) ,
\end{align*}
as our reference space throughout the paper, because the shifted spline coefficients are finite for these tempered distributions, see \eqref{eq:dualityBound}. For now, we only note that each space used in Theorems \ref{thm:mainBesov} and \ref{thm:mainTriebel} is continuously embedded into $\cB$ by standard embeddings; see \cite[Sections 2.3.2 and 2.7.1]{triebel1983function}.

\begin{thm}\label{thm:mainBesov}
Let $\{\psi,\Psi\}$ be a Battle-Lemari\'e system of order $n\in\nN_0$. Suppose $\frac{1}{2(n+1)} < p \leqslant \infty$ and $0<q\leqslant \infty$. Further, take
\begin{align}
    -\frac{1}{p'} - n < s < n+1 . \label{eq:sBesovRange}
\end{align}
Then $B^s_{p,q}(\nR)$ consists of those $f\in\cB$, such that
\begin{align}
    \left( \sum_{j=-1}^\infty 2^{j\left( s-\frac{1}{p} \right)q} \left( \sum_{\mu\in\nZ} |\sfk_{j,\mu}(f)|^p \right)^{q/p} \right)^{1/q} < \infty . \label{eq:BesovSplineNorm}
\end{align}
We make the usual modifications if $p$ or $q$ is $\infty$. Furthermore, for these values of $p,q,s$, the quantity in \eqref{eq:BesovSplineNorm} gives an equivalent (quasi-) norm on $B^s_{p,q}(\nR)$.
\end{thm}

Analogously, we have a result for the Triebel-Lizorkin spaces.

\begin{thm}\label{thm:mainTriebel}
Let $\{\psi,\Psi\}$ be a Battle-Lemari\'e system of order $n\in\nN_0$. Suppose $\frac{1}{2(n+1)} < p,q < \infty$ and
\begin{align}
    \max\left\{ -\frac{1}{p'} , -\frac{1}{q'} \right\} - n < s < n+1 . \label{eq:sTriebelRange}
\end{align}
In this parameter range, the space $F^s_{p,q}(\nR)$ is the collection of those $f\in\cB$ satisfying
\begin{align}
    \left\| \left( \sum_{j=-1}^\infty 2^{jsq} \left| \sum_{\mu\in\nZ} \sfk_{j,\mu}(f) 1_{I_{j,\mu}} \right|^q \right)^{1/q} \right\|_{L_p} < \infty . \label{eq:TriebelSplineNorm}
\end{align}
Furthermore, for these values of $p,q,s$, the quantity in \eqref{eq:TriebelSplineNorm} gives an equivalent (quasi-) norm on $F^s_{p,q}(\nR)$. The theorem also holds in the case $q=\infty$ (modifying \eqref{eq:TriebelSplineNorm} in the usual way) so long as $\frac{1}{2(n+1)} <p\leqslant 1$.
\end{thm}

\begin{rem}
We must avoid the case $q=\infty$ when $1<p<\infty$ in Theorem \ref{thm:mainTriebel} for a somewhat superficial reason: the currently known result, Theorem \ref{thm:SplineCharac}, is not sufficient to give this case. This issue appears when we are analyzing the validity of \eqref{eq:ChuiWangTriebel}.
\end{rem}

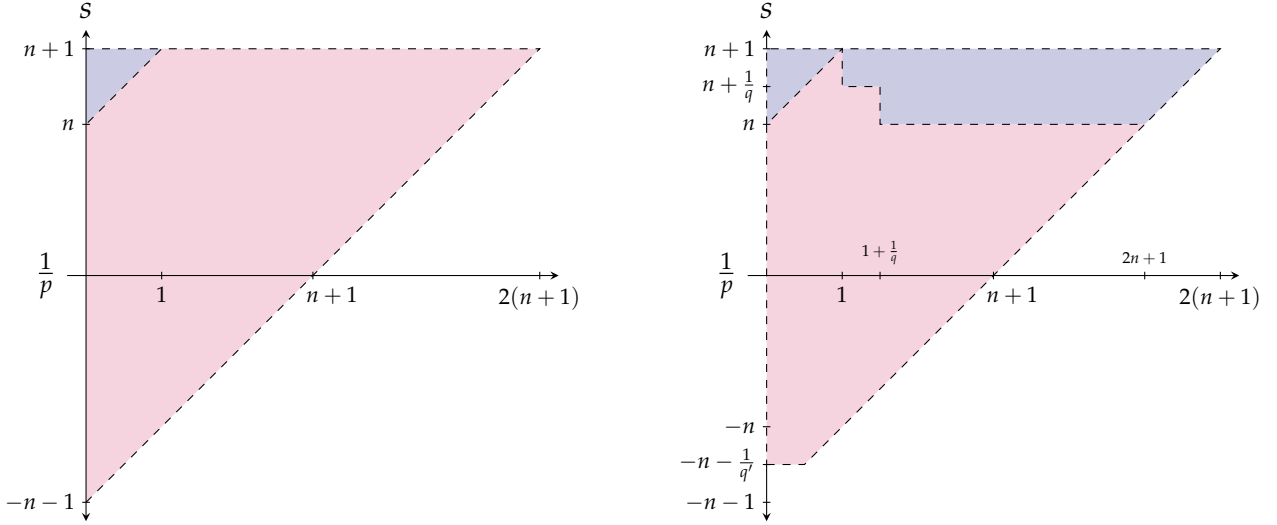
\begin{figure}[ht]
\centering

\begin{tikzpicture}[scale=1]

    \fill[fill=specialpink!70!white] (0,-3) -- (6,3) -- (1,3) -- (0,2) -- cycle; 
    \fill[fill=specialpurple!70!white] (0,2) -- (1,3) -- (0,3) -- cycle; 
    \draw[dashed] (0,-3) -- (6,3) -- (1,3) -- (0,2); 
    \draw[dashed] (0,3) -- (1,3); 

    \draw[-stealth] (-0.25,0) -- (6.25,0); 
    \node[anchor=east] at (-0.25,0) {$\frac{1}{p}$}; 
    \draw[stealth-stealth] (0,-3.25) -- (0,3.25); 
    \node[anchor=south] at (0,3.25) {$s$}; 
    \node[anchor=east] at (0,-3) {\scalebox{0.7}{$-n-1$}}; 
    \node[anchor=east] at (0,2) {\scalebox{0.7}{$n$}}; 
    \node[anchor=east] at (0,3) {\scalebox{0.7}{$n+1$}}; 
    \node[anchor=north] at (1,0) {\scalebox{0.7}{$1$}}; 
    \node[anchor=north] at (3.25,0) {\scalebox{0.7}{$n+1$}}; 
    \node[anchor=north] at (6,0) {\scalebox{0.7}{$2(n+1)$}}; 
    \draw (-0.05,-3) -- (0.05,-3); 
    \draw (-0.05,2) -- (0.05,2); 
    \draw (-0.05,3) -- (0.05,3); 
    \draw (1,-0.05) -- (1,0.05); 
    \draw (3,-0.05) -- (3,0.05); 
    \draw (6,-0.05) -- (6,0.05); 


    \fill[fill=specialpink!70!white] (9,-2.5) -- (9.5,-2.5) -- (14,2) -- (10.5,2) -- (10.5,2.5) -- (10,2.5) -- (10,3) -- (9,2) -- cycle; 
    \fill[fill=specialpurple!70!white] (9,2) -- (9,3) -- (10,3) -- cycle; 
    \fill[fill=specialpurple!70!white] (10,3) -- (15,3) -- (14,2) -- (10.5,2) -- (10.5,2.5) -- (10,2.5) -- cycle; 
    \draw[dashed] (9,-2.5) -- (9.5,-2.5) -- (14,2) -- (10.5,2) -- (10.5,2.5) -- (10,2.5) -- (10,3) -- (9,2) -- cycle; 
    \draw[dashed] (9,2) -- (9,3) -- (15,3) -- (14,2); 
    
    \draw[-stealth] (8.75,0) -- (15.25,0); 
    \node[anchor=east] at (8.75,0) {$\frac{1}{p}$}; 
    \draw[stealth-] (9,-3.25) -- (9,-2.5); 
    \draw[-stealth] (9,3) -- (9,3.25); 
    \node[anchor=south] at (9,3.25) {$s$}; 
    \node[anchor=east] at (9,-3) {\scalebox{0.7}{$-n-1$}}; 
    \node[anchor=east] at (9,2) {\scalebox{0.7}{$n$}}; 
    \node[anchor=east] at (9,3) {\scalebox{0.7}{$n+1$}}; 
    \node[anchor=east] at (9,-2) {\scalebox{0.7}{$-n$}};
    \node[anchor=east] at (9,2.5) {\scalebox{0.7}{$n+\frac{1}{q}$}};
    \node[anchor=east] at (9,-2.5) {\scalebox{0.7}{$-n-\frac{1}{q'}$}};
    \node[anchor=north] at (10,0) {\scalebox{0.7}{$1$}}; 
    \node[anchor=north] at (12.25,0) {\scalebox{0.7}{$n+1$}}; 
    \node[anchor=north] at (15,0) {\scalebox{0.7}{$2(n+1)$}}; 
    \node[anchor=south] at (10.5,0) {\scalebox{0.5}{$1+\frac{1}{q}$}}; 
    \node[anchor=south] at (14,0) {\scalebox{0.5}{$2n+1$}}; 
    \draw (8.95,-3) -- (9.05,-3); 
    \draw (8.95,2) -- (9.05,2); 
    \draw (8.95,3) -- (9.05,3); 
    \draw (8.95,-2) -- (9.05,-2); 
    \draw (8.95,-2.5) -- (9.05,-2.5); 
    \draw (8.95,2.5) -- (9.05,2.5); 
    \draw (10,-0.05) -- (10,0.05); 
    \draw (12,-0.05) -- (12,0.05); 
    \draw (15,-0.05) -- (15,0.05); 
    \draw (10.5,-0.05) -- (10.5,0.05); 
    \draw (14,-0.05) -- (14,0.05); 

\end{tikzpicture}

\caption{The parameter ranges for the Battle-Lemari\'e spline system to be an unconditional basis of $B^s_{p,q}$ (left) and $F^s_{p,q}$ (right, when $1<q<\infty$) are in pink, see Theorem \ref{thm:SplineCharac}. The extension provided by the frame in Theorems \ref{thm:mainBesov} and \ref{thm:mainTriebel} is in purple.} \label{fig:SplineFrame}
\end{figure}

We can pictorially compare the range for $s$ to be an unconditional basis, as given in \eqref{eq:sBesovIntro} and \eqref{eq:sTriebelIntro} (stated fully in Theorem \ref{thm:SplineCharac}), and the range, into which we have expanded using the spline frame. See Figure \ref{fig:SplineFrame}. For the $F$-spaces, we have only drawn the case $1<q<\infty$; similar images can easily be made for $0<q\leqslant 1$ and $q=\infty$ using Theorem \ref{thm:SplineCharac}.

As can be seen from Figure \ref{fig:SplineFrame}, both Theorems \ref{thm:mainBesov} and \ref{thm:mainTriebel} are only new for the larger smoothness values. The novelty here is a new description of the wavelet coefficients and the expansion of the smoothness range into $s$ values, where the Battle-Lemari\'e system does not give an unconditional basis. It is also a novelty in the area of wavelet characterizations for function spaces that the upper bound for $s$, i.e. $s<n+1$, has no dependence on $p$ or $q$.

We remark that it is difficult to exactly calculate spline coefficients involving the Battle-Lemari\'e system, such as those in Theorems \ref{thm:mainBesov} and \ref{thm:mainTriebel}, due to the system lacking compact support and not having an explicit representation (except in Fourier space). However, estimating these coefficients (e.g. in Lemmata \ref{lem:scale1} and \ref{lem:scale2}) becomes easier in many cases, because of the vanishing moments (property (IV) in Section \ref{subsec:wavelets}) commonly found in orthonormal wavelet systems. Furthermore, once we have proven results such as norm characterizations for the Battle-Lemari\'e system, we can use these results to easily prove analogous theorems for other wavelet systems coming from the same multiresolution analysis, e.g. the Chui-Wang wavelets \cite{chui1992compactly}. While these Chui-Wang wavelets are compactly supported, they do not give an orthogonal system and they lack vanishing moments, making our analysis in Sections \ref{sec:proofs} and \ref{sec:edgeCase} impossible. We will see in a future paper, however, how Theorems \ref{thm:mainBesov} and \ref{thm:mainTriebel} imply analogous results for the Chui-Wang system.

We conclude this introduction with one final question: can the ranges \eqref{eq:sBesovRange} and \eqref{eq:sTriebelRange} be expanded? In some cases we can easily show that the lower bound is optimal; we will briefly discuss this in Section \ref{sec:lowerOptimal}. The upper bound is more complicated, but we do have two results, which we discuss presently.

Below, $W^k_p(\nR)$ are the usual Sobolev spaces for $k\in\nN_0$ and $1\leqslant p\leqslant \infty$. For $k\in\nN$, we will use the norm
\begin{align}
    \|f\|_{W^k_p(\nR)} = \|f\|_{L_p(\nR)} + \|f^{(k)}\|_{L_p(\nR)} . \label{eq:sobolevNormDef}
\end{align}
Additionally, we will need to take a different reference space than $\cB$ for some technical reasons. We define $\cS_r\supset \cS$ for $r\in\nN_0$ to be the space of real-valued $\varphi\in C^r(\nR)$ satisfying for some $C_{\ell, k} > 0$
\begin{align*}
    |\varphi^{(k)}(x)| \leqslant C_{\ell, k} (1+|x|)^{-\ell}
\end{align*}
for $k=0,\ldots,r$, $\ell\in\nZ$, and $x\in\nR$. Convergence in $\cS_r$ is defined in a way analogous to that in $\cS$. We write $\cS_r'$ to mean the topological dual of $\cS_r$. Note that if $\{\psi,\Psi\}$ is a Battle-Lemari\'e system of order $n\geqslant 1$, then $\psi,\Psi\in \cS_{n-1}$ by properties (I) and (III) in Section \ref{subsec:wavelets}.


\begin{thm}\label{thm:edgeCase}
Suppose $n\geqslant 1$ and $\{\psi,\Psi\}$ is a Battle-Lemari\'e system of order $n$. Let $1 < p \leqslant \infty$. Then $W^{n+1}_p(\nR)$ consists of those $f\in \cS_{n-1}'(\nR)$, such that
\begin{align}
    \sup_{j\geqslant -1} 2^{j\left(n+1 - \frac{1}{p}\right)} \left( \sum_{\mu\in\nZ} |\sfk_{j,\mu}(f)|^p \right)^{1/p} < \infty ,\label{eq:SobolevSplineNorm}
\end{align}
making the usual modification if $p=\infty$. The above quantity is equivalent to the norm on $W^{n+1}_p(\nR)$ given in \eqref{eq:sobolevNormDef}.
\end{thm}

One might think that the proper way to recreate the Sobolev norm would be to use the sequence space norm $f^{n+1}_{p,2}$ (given in Definition \ref{def:sequenceSpace}), since $W^{n+1}_p = F^{n+1}_{p,2}$ when $1<p<\infty$, but above we have used $b^{n+1}_{p,\infty}$. This actually has precedence, however, as noted in \cite{garrigos2023haar}, results of this type have appeared in the literature as early as the 1960s, see, e.g., \cite{bovckarev1969coefficients}.

We will actually show slightly more than what the above theorem claims. In particular, we will see that some embeddings still hold when $p=1$. For $1\leqslant p\leqslant \infty$ and $n\geqslant 1$ let 
\begin{align*}
    B^{n+1,BL}_{p,\infty}(\nR) = \left\{ f\in \cS_{n-1}' : \eqref{eq:SobolevSplineNorm} \text{ holds} \right\} ,
\end{align*}
adjusting \eqref{eq:SobolevSplineNorm} if $p=\infty$. Then, we have the following.

\begin{prop}\label{prop:p1inftyEmbed}
If $1\leqslant p\leqslant\infty, n\geqslant 1$, then $W^{n+1}_p(\nR)\hookrightarrow B^{n+1,BL}_{p,\infty}(\nR)$ is a continuous embedding.
\end{prop}

In \cite[Theorem 1.4]{garrigos2023haar} the case $n=0$ was already proven for Theorem \ref{thm:edgeCase} and the authors were able to include a characterization of $BV$, the space of functions of bounded variation, when $p=1$. We will make some comments about $p=1$ and characterizing $BV$ in Remark \ref{rem:BV}. Additionally, in \cite[Theorem 1.4]{garrigos2023haar}, instead of taking $\cS_{n-1}'$ as the reference space as we did in Theorem \ref{thm:edgeCase}, the authors continued with $\cB$. They were able to do this, because they had the benefit of a previously developed, fairly robust, theory of some projections corresponding to their multiresolution analysis.\footnote{We define these projections in \eqref{eq:projection} and their version is equation (3.4) in \cite{garrigos2023haar}.} Developing this theory in our general setting would be a significant departure from the purpose of this paper, so we accept this small difference.

There are a number of other obvious questions to ask, such as what happens when $s > n+1$. One result in this direction, \cite[Proposition 11.1]{garrigos2023haar}, is that when $n=0$ and $s > 1$ or $s=1$ and $q=\infty$, the only continuously differentiable functions with \eqref{eq:BesovSplineNorm} finite are constant functions. Such a result for general $n\in\nN$ is currently not available.

\begin{rem}
Concerning the optimality of the ranges for $p$ and $q$ in Theorems \ref{thm:mainBesov} and \ref{thm:mainTriebel}, the only reason we could not take $0 < p,q\leqslant\infty$ ($p<\infty$ in the Triebel-Lizorkin case) is because \eqref{eq:sBesovRange} and \eqref{eq:sTriebelRange} are empty outside of the given ranges for $p,q$. In particular, if one somehow managed to expand \eqref{eq:sBesovRange} and \eqref{eq:sTriebelRange}, the ranges for $p,q$ would also expand for free.
\end{rem}


\noindent\textbf{Notation.} We write $X\lesssim Y$ to mean there is a constant $C > 0$ such that $X\leqslant CY$. By $X\sim Y$ we will mean that $X\lesssim Y$ and $Y\lesssim X$ hold simultaneously, i.e. there is a constant $C > 0$ so that $C^{-1}Y \leqslant X\leqslant CY$. This constant may depend on parameters, e.g. $p$, in which case we will sometimes write $X\lesssim_p Y$.

A few unorganized pieces of notation are as follows: $1_S$ means the indicator function on a measurable set, $S$, the H\"older conjugate (or dual) of an integration parameter, say $p$, will be written as $p'$, i.e. it satisfies $\frac{1}{p} + \frac{1}{p'} = 1$ (if $p=1$, then we define $p'=\infty$),\footnote{Sometimes in the literature one defines $p' = \infty$ if $0<p<1$. We do \emph{not} do that here.} $\delta_{j,k}$ is the Kronecker-$\delta$ function, and $j_+ = \max\{0,j\}$ for any real $j$. As usual, $L_p(\nR)$ for $0<p\leqslant\infty$ will mean the space of $p$-integrable, measurable functions and $L_1^{\text{loc}}(\nR)$ is the space of locally integrable, measurable functions. Also, $\cS$ is the Schwartz space and $\cS'$ is its dual, the space of tempered distributions. A superscript with parentheses will be used to denote derivatives, e.g. $f^{(k)}$ indicates the $k$-th derivative of $f$ (it will be clear from context if these derivatives are classical or distributional). In accordance with what is fairly standard in wavelet analysis, we choose the following Fourier transform:
\begin{align*}
    \cF f(\xi) = \widehat{f}(\xi) = \int_\nR f(x) e^{- ix\xi} \, dx .
\end{align*}
The inverse Fourier transform will be notated as $\cF^{-1}$. We use the notation $f*g$ for convolution, for $f,g$ such that the convolution is defined. More specialized notation is in Section \ref{sec:prereq}.
\vspace{+1mm}

\noindent\textbf{Idea of the Proof.} We shall see that the skeleton of the proof for our theorems is completely analogous to that given in \cite{garrigos2023haar} when $n=0$; what distinguishes this paper is the details. There are two main hurdles, which complicate our proofs: the Battle-Lemari\'e system has no compact support (only exponential decay) and it has no closed form (in physical space) that we can write down and manipulate. In Section \ref{sec:weak}, this is particularly apparent when one compares \cite[Lemma 4.2]{garrigos2023haar} to our versions: Lemmata \ref{lem:scale1}, \ref{lem:scale2convolution}, and \ref{lem:scale2}. Moving to Section \ref{sec:strong}, the striking simplicity of \cite[Lemma 2.1]{garrigos2023haar} is contrasted by the excursion we must take through \cite{ushakova2018localisation} in Section \ref{subsubsec:remarkableConnection} to prove Lemma \ref{lem:hardPartCrux}. These roadblocks pose significant issues in Section \ref{sec:edgeCase} as well, the most apparent manifestation of which is the need to prove Lemma \ref{lem:rhoExpDecay}, which is perfectly obvious when $n=0$.

\vspace{+1mm}

\noindent\textbf{Plan of the Paper.} We begin in Section \ref{sec:prereq} with a review of the various definitions and results that we will need from the literature. This will include comments on splines and wavelets, basic definitions and properties of Besov and Triebel-Lizorkin spaces, the relation of spline wavelets to these spaces, some results on maximal functions, and a clarification of the use of the space $\cB$ in our main theorems. Following this, in Section \ref{sec:proofs}, we prove Theorems \ref{thm:mainBesov} and \ref{thm:mainTriebel} in tandem, starting with the weaker direction in Section \ref{sec:weak} and later moving to the stronger in Section \ref{sec:strong}. We, then, conclude with a proof of Theorem \ref{thm:edgeCase} and Proposition \ref{prop:p1inftyEmbed} in Section \ref{sec:edgeCase}.
\vspace{+1mm}

\noindent\textbf{Acknowledgments.} This paper was born from the author's master's thesis; accordingly, he would like to express his deep gratitude towards his two master's advisors: Christoph Thiele and Rajula Srivastava. Their advice on this paper and on math in general has been indispensable; further, the author would like to thank them for suggesting such an interesting problem. He would be remiss without further expressing his appreciation for the unending support that R. Srivastava has given even after the completion of his master's thesis---it has made all the difference. The author would also like to thank the two anonymous reviewers for their meticulous reading of the manuscript and their many helpful comments.

Additionally, the author was supported by a grant from the Deutscher Akademischer Austauschdienst (DAAD) throughout his master's degree; this work would not have been possible without their generosity. During the writing of this paper, the author was supported by the Primus research programme PRIMUS/21/SCI/002 of Charles University. This research was supported in part by the grant no. 23-04720S of the Czech Science Foundation. The author also acknowledges support from the Charles University Research Centre program UNCE/24/SCI/005. Finally, the author was supported by the grant SVV-2025-260827.

%% file: Sections/preliminaries.tex
\section{Prerequisite Material}\label{sec:prereq}


\subsection{Splines and Wavelets}\label{subsec:wavelets}

The multiresolution analysis underlying this paper corresponds to the Battle-Lemari\'e wavelet system, but its origin is in the B-splines. One can find the definition of a multiresolution analysis in \cite[Definition 2.2]{wojtaszczyk1997mathematical}. We will write $(V_j)_{j\in\nZ}$ to denote the scale spaces of a multiresolution analysis for $L_2(\nR)$ and $W_j$ will be the orthogonal complement of $V_j$ in $V_{j+1}$, i.e.
\begin{align*}
    V_{j+1} = V_j \oplus W_j , \quad j\in\nZ .
\end{align*}
We say that $G\in V_0$ is a scaling function if the system $(G(x-\mu))_{\mu\in\nZ}$ gives a Riesz basis of $V_0$; similarly, we will say $g\in W_0$ is a wavelet corresponding to $G$ if $(g(x-\mu))_{\mu\in\nZ}$ gives a Riesz basis of $W_0$. Throughout the paper, we will use the notation $V_j$ and $W_j$ for scale and wavelet spaces, respectively; it will be made clear each time to which scaling function/wavelet these correspond.

A classic example of a scaling function is, of course, the B-spline, $B_m$, of order $m\in\nN_0$. These are given by
\begin{align}
    B_0(x) = 1_{[0,1)}(x) , \quad B_m(x) = (B_{m-1} * 1_{[0,1)})(x) \label{eq:BsplineDef}
\end{align}
for $m\in\nN$.\footnote{There are competing notations in the literature. Some authors write $N_m$ instead of $B_m$ and some start $m$ at 1, i.e. $B_1 = 1_{[0,1)}$.} There are many interesting properties of B-splines, but the only two we will need are that $B_m$ has compact support for each $m\in\nN_0$ and the remarkable derivative identity
\begin{align}
    B_m'(x) = B_{m-1}(x) - B_{m-1}(x-1) , \label{eq:magicIdentity}
\end{align}
which holds for $m\in\nN$. For proofs of these facts, along with many other results, check \cite[Theorem 4.3]{chui1992introduction}. 

What the B-splines lack is intrascale orthogonality, i.e. if $j=k$, then $(B_{m;j,\mu} , B_{m;k,\nu})$ is not zero for all $\mu,\nu\in\nZ$ with $\mu\neq\nu$ (we define $B_{m;j,\mu}(x) = B_m(2^j x - \mu)$). This is what prompts us to orthonormalize, creating the Battle-Lemari\'e scaling function, $\Psi$, and, from that, its corresponding wavelet, $\psi$. Both of these are real-valued and they enjoy a number of important properties, which we list below for $\psi,\Psi$ of order $n\in\nN_0$.
\begin{enumerate}[label=(\Roman*)]
    \item Globally, $\psi,\Psi\in C^{n-1}(\nR)$ for $n\in\nN$. If $n=0$, they are measurable.
    \item On each interval $(\mu,\mu+1/2)$ with $\mu\in\frac12 \nZ$, $\psi$ and $\Psi$ are polynomials of degree at most $n$.
    \item There exist $C_0,\gamma > 0$ (which may depend on $n$), such that
    \begin{align*}
        |\Psi^{(k)}(x)| + |\psi^{(k)}(x)| \leqslant C_0 e^{-\gamma|x|}
    \end{align*}
    for $0\leqslant k\leqslant n-1$ and all $x\in\nR$. If $n=0$, we also require that $\psi,\Psi$ have compact support.
    \item The wavelet, $\psi$, has $n+1$ vanishing moments, i.e.
    \begin{align*}
        \int_\nR x^\kappa \psi(x) \, dx = 0
    \end{align*}
    for $\kappa = 0,1,\ldots,n$.
\end{enumerate}
We will often describe a function as having exponential decay; for example, we have done so without reference to this terminology for $\psi,\Psi$ in property (III) above. We define this explicitly.

\begin{define}\label{def:expDecay}
A function, $f : X \to \nF$, where $X$ is a subset of either $\nC$ or $\nZ$ and $\nF\in \{\nR,\nC\}$, is said to decay exponentially or have exponential decay if there exist constants, $C,c > 0$, such that
\begin{align*}
    |f(x)| \leqslant C e^{-c |x|}
\end{align*}
for all $x\in X$.
\end{define}

While there is a standard function, which is the Battle-Lemari\'e scaling function, see \eqref{eq:classicalBLScaling}, and likewise for the wavelet, see \eqref{eq:classicalBLWavelet}, there are many choices that need to be made along the path of orthonormalization, see \cite{ushakova2018localisation}; we handle this with the following definition.

\begin{define}[Battle-Lemari\'e Wavelet System]\label{def:bl}
A pair of functions, $\{g,G\}$, consisting of a scaling function, $G$, and a wavelet, $g$, will be called a Battle-Lemari\'e wavelet system of order $m\in\nN_0$ if the following conditions hold. Write $g_{j,\mu}(x) = g(2^j x - \mu)$ for $j\in\nN_0$, $\mu\in\nZ$ and $g_{-1,\mu}(x) = \sqrt{2} G(x-\mu)$ for $\mu\in\nZ$.
\begin{enumerate}[label=(\roman*)]
    \item The system has both inter- and intrascale orthonormality, i.e.
    \begin{align}
        (g_{j,\mu},g_{k,\nu}) = 2^{-j} \delta_{j,k} \delta_{\mu,\nu} \label{eq:orthodef}
    \end{align}
    for $j,k\in\nN_0\cup\{-1\}$ and $\mu,\nu\in\nZ$.
    \item  It must satisfy properties (I)-(IV) above.
    \item It must correspond to the same multiresolution analysis of $L_2$ as $B_m$, i.e. $\{G(\cdot - \mu) : \mu\in\nZ\}$ is an orthonormal basis of the scale space, $V_0$, from the scaling function $B_m$, and $\{g(\cdot-\mu) : \mu\in\nZ\}$ is an orthonormal basis of the wavelet space, $W_0$, corresponding to this $V_0$.
\end{enumerate}
\end{define}

\begin{rem}\label{rem:Haar}
When $n=0$, the classical Battle-Lemari\'e system of order zero is the Haar system, $\{h,H\}$, which was first investigated by A. Haar\footnote{This is Alfr\'ed Haar, not the author of this paper, Andrew Haar.} in the early 20th century \cite{haar1909theorie}. They are defined as
\begin{align}
    h(x) = 1_{[0,1/2)}(x) - 1_{[1/2,1)}(x) , \quad H(x) = 1_{[0,1)}(x) . \label{eq:HaarSystem}
\end{align}
It will be useful to note that, by \cite[Lemma 4.4]{wojtaszczyk1997mathematical}, two compactly supported functions, which both give rise to an orthonormal basis of the same subspace of $L_2(\nR)$ through their integer shifts, must themselves be integer shifts of each other. That is to say, if $\{\psi,\Psi\}$ is a Battle-Lemari\'e system of order zero, then there exist $k,\ell\in\nZ$ such that $\psi(x) = h(x-k)$ and $\Psi(x) = H(x-\ell)$.
\end{rem}

\begin{rem}\label{rem:classicalSystem}
We define the classical Battle-Lemari\'e system for $n\in\nN$. As in \cite[Example 7.1]{chui1992introduction} or \cite[Chapter 5]{daubechies1992ten}, if we temporarily say $\{\psi,\Psi\}$ is the classical system of order $n$,
\begin{align}
    \widehat{\Psi}(\xi) = \left( \sum_{k\in\nZ} |\widehat{B}_n(\xi + 2\pi k)|^2 \right)^{-1/2} \widehat{B}_n(\xi) , \label{eq:classicalBLScaling}
\end{align}
with the series converging pointwise, and 
\begin{align}
    \widehat{\psi}(\xi) = - e^{-i\xi/2} \frac{\overline{\widehat{\Psi}(\xi + 2\pi)}}{\overline{\widehat{\Psi}(\xi/2 + \pi)}} \widehat{\Psi}(\xi/2) . \label{eq:classicalBLWavelet}
\end{align}
The system above gives the Haar system from Remark \ref{rem:Haar} when $n=0$, see \cite[Example 7.1]{chui1992introduction}.
\end{rem}

The following lemma captures the interplay of properties (I)-(III) of our Battle-Lemari\'e system. Indeed, because we require the polynomials defining $\psi,\Psi$ to connect in such a way that they are together in $C^{n-1}(\nR)$, the non-leading coefficients of adjacent polynomial pieces when written around their center are equal. The lemma was first proven in \cite{srivastava2023orthogonal} in the process of proving the necessity of the smoothness range \eqref{eq:sTriebelIntro} for the Battle-Lemari\'e system to be an unconditional basis of $F^s_{p,q}(\nR)$.

\begin{lem}[Lemma 3.1 from \cite{srivastava2023orthogonal}]\label{lem:psiPoly}
Let $\{\psi,\Psi\}$ be a Battle-Lemari\'e system of order $n\in\nN_0$ and fix $\mu\in\nZ$. Say
\begin{align*}
    \psi(x) = A_{\mu-1}^n \left( x - \frac{\mu}{2} \right)^n + A^{n-1}_{\mu-1} \left( x - \frac{\mu}{2} \right)^{n-1} + \cdots + A_{\mu-1}^0
\end{align*}
on $[\frac{\mu - 1}{2}, \frac{\mu}{2}]$ and
\begin{align*}
    \psi(x) = A_\mu^n \left( x - \frac{\mu}{2} \right)^n + A^{n-1}_\mu \left( x - \frac{\mu}{2} \right)^{n-1} + \cdots + A_\mu^0
\end{align*}
on $[\frac{\mu}{2},\frac{\mu+1}{2}]$. Then $A_{\mu-1}^k = A_\mu^k$ for each $k = 0,1,\ldots,n-1$ and
\begin{align}
    |A^n_\mu| \leqslant 4C_0 e^{\gamma/2} e^{-(\gamma/2) |\mu|} , \quad \mu\in\nZ, \label{eq:leadingCoefBound}
\end{align}
where the $C_0,\gamma$ above are the same as in property (III) of $\psi$. The entire lemma also holds for $\Psi$.
\end{lem}

\begin{rem}
The lemma was proven in \cite{srivastava2023orthogonal} for the classical Battle-Lemari\'e wavelet system, but if one examines the proof, the only properties of $\psi,\Psi$ that are used are (I)-(III) above.
\end{rem}


\subsection{Multiresolution Analyses for Spaces of Distributions}\label{sec:multiresDistribution}

The following is a summary of \cite[Section 3]{walter1992wavelets}. Recall the definition of $\cS_r$ from the Introduction. We fix the Battle-Lemari\'e system $\{\psi,\Psi\}$ of order $n\geqslant 1$ for this subsection.

Because both $\Psi$ and $\psi$ are in $\cS_{n-1}$, the wavelet coefficients from the base scales, i.e. testing against $(\Psi(\cdot-\mu))$ or $(\psi(\cdot-\mu))$, for distributions in $\cS_{n-1}'$ will be $O(|\mu|^k)$ for some $k\in\nZ$. The expansion of some $f\in\cS_{n-1}'$ at this base scale, then, makes sense and converges in the sense of $\cS_{n-1}'$.

\begin{define}
Say $(a_\mu)\subset\nC$ is taken to be such that $a_\mu = O(|\mu|^k)$ for some $k\in\nZ$ and define
\begin{align*}
    T_0 = \left\{f : f(x) = \sum_{\mu\in\nZ} a_\mu \Psi(x-\mu)\right\} , \quad U_0 = \left\{f : f(x) = \sum_{\mu\in\nZ} a_\mu \psi(x-\mu) \right\} .
\end{align*}
\end{define}

The corresponding dilation spaces are denoted with $T_m$ and $U_m$. With these preparations, we have a multiresolution analysis of $\cS_{n-1}'$. In particular, following \cite{walter1992wavelets},
\begin{align*}
    \{0\}\subset \cdots \subset T_{-1}\subset T_0\subset T_1\subset \cdots \subset \cS_{n-1}'
\end{align*}
and
\begin{align*}
    \bigcap_{m\in\nZ} T_m = \{0\} , \quad \overline{\bigcup_{m\in\nZ} T_m} = \cS_{n-1}' ,
\end{align*}
with the closure being in the topology of $\cS_{n-1}'$. Further, it holds
\begin{align}
    T_m = T_0 \dotplus U_0 \dotplus U_1\dotplus \cdots \dotplus U_{m-1} , \label{eq:mysavior}
\end{align}
where $\dotplus$ indicates a direct (but not orthogonal) sum.\footnote{In \cite[Theorem 5]{walter1992wavelets} the result is written as $T_m = T_0 \dotplus U_0\dotplus \cdots \dotplus U_m$, but this is a typo. The line before says ``$U_m$ are the complementary subspaces of $T_m$ in $T_{m+1}$;'' in other words, $T_{m+1} = T_m \dotplus U_m$, i.e. \eqref{eq:mysavior} is the correct expression.} Furthermore, if $E_m : \cS_{n-1}'\to T_m$ is the projection from $\cS_{n-1}'$ onto $T_m$, then $E_m$ is continuous and $E_m f \to f$ in $\cS_{n-1}'$ for $f\in\cS_{n-1}'$.


\subsection{Besov and Triebel-Lizorkin Spaces}\label{subsec:spaces}

This subsection will cover the basic facts about the Besov and Triebel-Lizorkin spaces, written $B^s_{p,q}(\nR)$ and $F^s_{p,q}(\nR)$, respectively. Both are Banach spaces when $1\leqslant p,q\leqslant \infty$ (and $s\in\nR$) and quasi-Banach spaces otherwise.

We give a definition of the Besov and Triebel-Lizorkin spaces through Littlewood-Paley pieces with vanishing moments. Consider pieces $\beta_0, \beta\in \cS$ such that, for a fixed $\varepsilon > 0$,
\begin{align}
    \begin{cases}
        |\widehat{\beta}_0(\xi)| > 0 & |\xi| < 2\varepsilon , \\
        |\widehat{\beta}(\xi)| > 0 & \frac{\varepsilon}{2} < |\xi| < 2\varepsilon ,
    \end{cases} \label{eq:betaNonzero}
\end{align}
and write $\beta_k(x) = 2^k \beta(2^k x)$ if $k\geqslant 1$. We also require that $\beta$ has $m\in\nN_0$ vanishing moments, i.e. $\widehat{\beta}^{(\kappa)}(0) = 0$ for $0\leqslant \kappa \leqslant m-1$ (if $m=0$, then the condition is empty). Define
\begin{align*}
    L_k f = \beta_k * f
\end{align*}
for $k\in\nN_0$ and $f\in\cS'(\nR)$.

We may now write the $B^s_{p,q}$ and $F^s_{p,q}$ quasi-norms in terms of those for $\ell_q(L_p)$ and $L_p(\ell_q)$ for $0<p,q\leqslant\infty$, which we define here for completeness. Each is a space of sequences of complex-valued, Lebesgue measurable functions on $\nR$, say $(g_j(x))_j$, with the following (quasi-) norms:
\begin{align*}
    \|(g_j)\|_{L_p(\ell_q)} &= \left( \int_\nR \left( \sum_j |g_j(x)|^q \right)^{p/q} \, dx \right)^{1/p} , \\
    \|(g_j)\|_{\ell_q(L_p)} &= \left( \sum_j \left( \int_\nR |g_j(x)|^p \, dx \right)^{q/p} \right)^{1/q} ,
\end{align*}
making the usual modifications if $p$ or $q$ is infinity. We recall that $L_p, \ell_p, L_p(\ell_q), \ell_q(L_p)$ are all $u$-norms; in the former two, one takes $u = \min\{p,1\}$, while $u = \min\{p,q,1\}$ in the latter two.

We have the following definitions. Given $s<m$ and $0<p,q\leqslant\infty$, the Besov space, $B^s_{p,q}(\nR)$, is those $f\in\cS'(\nR)$, such that
\begin{align}
    \|f\|_{B^s_{p,q}(\nR)} := \|(2^{ks} L_k f)_{k\in\nN_0}\|_{\ell_q(L_p)} < \infty . \label{eq:normDefBesov}
\end{align}
For $s < m$, $0<q\leqslant\infty$, $0<p<\infty$, the Triebel-Lizorkin space, $F^s_{p,q}(\nR)$, consists of the tempered distributions, $f\in\cS'(\nR)$, such that
\begin{align}
    \|f\|_{F^s_{p,q}(\nR)} := \|(2^{ks} L_k f)_{k\in\nN_0}\|_{L_p(\ell_q)} < \infty . \label{eq:normDefTriebel}
\end{align}
There are many equivalent definitions that can be used for these spaces, see \cite{triebel1983function,triebel2006function}.

\begin{rem}\label{rem:equivDef}
Each choice of $\beta_0,\beta$ gives rise to an equivalent (quasi-) norm on these spaces \cite{rychkov1999theorem}.
\end{rem}

We will take $\varepsilon = 1/2$ in \eqref{eq:betaNonzero} and the number of vanishing moments of $\beta$ to be $n+1$. Furthermore, we may assume that
\begin{align*}
    \supp\beta_0,\supp\beta \subset (-1/4,1/4) .
\end{align*}
This will be useful in Section \ref{sec:weak}. 

It would also be useful to have some kind of decomposition of a distribution in terms of these Littlewood-Paley pieces. To mend this absence, we require slightly more of $\beta_0,\beta$: we will assume that $|\widehat{\beta}_0(\xi)| > 0$ for $|\xi|\leqslant 1$ and $|\widehat{\beta}(\xi)| > 0$ for $\frac14 \leqslant |\xi| \leqslant 1$. Now, let $\eta$ be a smooth function with $\eta(\xi) = 1$ on $\{|\xi| \leqslant 1/2\}$, $|\eta(\xi)| > 0$ on $\{|\xi| < 1\}$ and $\supp\eta = \{|\xi|\leqslant 1\}$. Write
\begin{align}
    \widehat{\Lambda_0}(\xi) = \frac{\eta(\xi)}{\widehat{\beta}_0(\xi)} \quad \text{and} \quad \widehat{\Lambda}(\xi) = \frac{\eta(\xi) - \eta(2\xi)}{\widehat{\beta}(\xi)} . \label{eq:LambdaDef}
\end{align}
As usual, $\widehat{\Lambda_k}(\xi) = \widehat{\Lambda}(2^{-k}\xi)$ for $k\geqslant 1$ and $\Lambda_k f = \Lambda_k * f$ for $k\geqslant 0$ and $f\in\cS'$. Our decomposition identity is, then,
\begin{align}
    f = \sum_{k=0}^\infty L_k (\Lambda_k f) \label{eq:LPdecomp}
\end{align}
for all $f\in\cS'(\nR)$ with convergence in $\cS'(\nR)$. This decomposition has appeared in the literature many times before, e.g. \cite[Equation (2.5)]{garrigos2023haar} or \cite[Equation (10)]{rychkov1999theorem}. See also \cite[Lemma 3.6]{ullrich2016role}.

We have, thus, given ourselves some flexibility as well by breaking this decomposition of $f$ into parts: $L_k$ and $\Lambda_k$. Conveniently,
\begin{align}
    \|(2^{ks}\Lambda_k f)\|_{\ell_q(L_p)} \sim_{\Lambda,\Lambda_0} \|f\|_{B^s_{p,q}} , \quad \|(2^{ks} \Lambda_k f)\|_{L_p(\ell_q)} \sim_{\Lambda,\Lambda_0} \|f\|_{F^s_{p,q}} \label{eq:LambaBound}
\end{align}
for $s\in\nR$ and $0<p,q\leqslant \infty$ (requiring $p<\infty$ for $F^s_{p,q}$). This is because $\Lambda$ satisfies the properties in \eqref{eq:betaNonzero} with $\varepsilon = 1/2$ and has infinitely many vanishing moments (also recall Remark \ref{rem:equivDef}).

Coming back to our Battle-Lemari\'e system, there is a characterization of the Besov and Triebel-Lizorkin spaces through these wavelets. Naturally, the idea is to give an isomorphism through wavelet coefficients to some sequence space, which we define presently.

\begin{define}[Sequence Spaces]\label{def:sequenceSpace}
Let $0<p,q\leqslant \infty$, $s\in\nR$, and take $\omega = (\omega_{j,\mu})_{j\geqslant -1,\mu\in\nZ}\subset \nC$. On the space of such sequences, we have the following (quasi-) norms:
\begin{align*}
    \|\omega\|_{b^s_{p,q}} &= \left( \sum_{j=-1}^\infty \left( 2^{j\left( s - \frac{1}{p} \right)} \left( \sum_{\mu\in\nZ} |\omega_{j,\mu}|^p \right)^{1/p} \right)^q \right)^{1/q} , \\
    \|\omega\|_{f^s_{p,q}} &= \left\| \left( \sum_{j=-1}^\infty 2^{jsq} \left| \sum_{\mu\in\nZ} \omega_{j,\mu} 1_{I_{j,\mu}}(\cdot) \right|^q \right)^{1/q} \right\|_{L_p} ,
\end{align*}
requiring $p<\infty$ for $f^s_{p,q}$. If $p$ or $q$ is infinity, one makes the usual modifications.
\end{define}

The theorem below comes from \cite[Theorems 2.46 and 2.49]{triebel2010bases}.\footnote{The sequence spaces in \cite{triebel2010bases} are $b_{p,q}^-$ and $f_{p,q}^-$ (equations (2.100)-(2.102)), but these are related to Definition \ref{def:sequenceSpace} in an obvious way.}

\begin{thm}\label{thm:SplineCharac} Let $\{\psi,\Psi\}$ be a Battle-Lemari\'e system of order $n\in\nN_0$.
\begin{enumerate}[label=(\roman*)]
    \item Let $0<p,q\leqslant\infty$ and
    \begin{align}
        -\frac{1}{p'} - n < s < n + \min \left\{ \frac{1}{p} , 1 \right\} . \label{eq:sOGBesovRange}
    \end{align}
    Then $f\in \cS'$ is in $B^s_{p,q}(\nR)$ if and only if it can be represented as 
    \begin{align}
        f = \sum_{j = -1}^\infty \sum_{\mu\in\nZ} \omega_{j,\mu} \psi_{j,\mu} \label{eq:SplineSeriesRep}
    \end{align}
    for some $\omega\in b^s_{p,q}$. The sequence, $\omega$, is uniquely determined by
    \begin{align}
        \omega_{j,\mu}(f) = 2^j (f,\psi_{j,\mu}) . \label{eq:coeffIsoSpline}
    \end{align}
    The unconditional convergence in \eqref{eq:SplineSeriesRep} takes place in $\cS'(\nR)$. The map $f\mapsto \omega(f)$ gives an isomorphism from $B^s_{p,q}$ to $b^s_{p,q}$.
    \item Let $0<p<\infty$, $0<q\leqslant\infty$, and
    \begin{align}
        \max\left\{ -\frac{1}{p'}, -\frac{1}{q'} , 0 \right\} - n < s < n . \label{eq:sOGTriebelRange}
    \end{align}
    Then $f\in\cS'$ is in $F^s_{p,q}$ if and only if it can be written as in \eqref{eq:SplineSeriesRep} with $\omega\in f^s_{p,q}$ uniquely given by \eqref{eq:coeffIsoSpline}. Convergence for \eqref{eq:SplineSeriesRep} is in $\cS'$ and $f\mapsto \omega(f)$ is an isomorphism from $F^s_{p,q}$ to $f^s_{p,q}$.
    \item If, instead, $1<p,q<\infty$, then item (ii) holds for the improved range
    \begin{align}
        \max\left\{ -\frac{1}{p'} , -\frac{1}{q'} \right\} - n < s < n + \min\left\{ \frac{1}{p} , \frac{1}{q} \right\} . \label{eq:sTriebel>1}
    \end{align}
    This also holds if $1\leqslant \frac{1}{p} < 1 + \frac{1}{q}$ given $0<p\leqslant 1$ and $1<q<\infty$.
\end{enumerate}
\end{thm}


\begin{rem}
If $0 < p \leqslant 1$, then Theorem \ref{thm:SplineCharac} actually gives us Theorem \ref{thm:mainBesov} without the oversampling. Of course, if $0 < p \leqslant \frac{1}{2(n+1)}$, then \eqref{eq:sOGBesovRange} is empty, so it does not give us any extra result.
\end{rem}


\subsection{Maximal Functions}\label{subsec:maximal}

The most classic maximal function is the Hardy-Littlewood maximal function, which, on $\nR$ is given for $g\in L_1^{\text{loc}}(\nR)$ by
\begin{align*}
    Mg(x) = \sup_{I\ni x} \frac{1}{|I|} \int_I |g(y)| \, dy ,
\end{align*}
the supremum being over (finite) intervals, which contain $x$. 

The Fefferman-Stein inequality, first proven in \cite{fefferman1971some}, says that if $1 < p < \infty$ and $1 < q\leqslant \infty$, then $\|(Mg_j)\|_{L_p(\ell_q)}\lesssim_{p,q} \|(g_j)\|_{L_p(\ell_q)}$, where $(g_j) \in L_p(\ell_q)$. We can lift the requirement that $p,q>1$ with the following: define for $g\in L_1^{\text{loc}}(\nR)$ and $r > 0$
\begin{align*}
    M_r g(x) = \left( M(|g|^r)(x) \right)^{1/r} .
\end{align*}
A simple corollary of the Fefferman-Stein inequality is that if $0<p<\infty$, $0<q\leqslant\infty$, and $0<r<\min\{p,q\}$, then for any sequence, $(g_j)\in L_p(\ell_q)$,
\begin{align}
    \|(M_r g_j)\|_{L_p(\ell_q)} \lesssim_{p,q,r} \|(g_j)\|_{L_p(\ell_q)} . \label{eq:BetterFeffermanStein}
\end{align}

We will now use this to prove a useful result for later. Recall the definitions of the shifted spline coefficients from \eqref{eq:shiftedSplineCoef} and \eqref{eq:shiftedSplinej=-1} and the dyadic intervals, $I_{j,\mu}$, from \eqref{eq:dyadicInterval}. The following lemma is not new, see \cite[Section 1.5.3]{triebel2006function}, but we could not find a proof in the literature. Our method of proof takes ideas from \cite[Lemma 4.6]{garrigos2023haar}.

\begin{lem}\label{lem:hardDirecMaxEstim}
Fix $m\in\nZ$, $0<p<\infty$, $0<q\leqslant\infty$, and $s\in\nR$. If $0 < r < \min\{p,q\}$, then for all $f\in\cB$,
\begin{align}
    \left\| \left( \sum_{j=-1}^\infty 2^{jsq} \left| \sum_{\mu\in\nZ} \sfk_{j,\mu+m}(f) 1_{I_{j,\mu}} \right|^q \right)^{\frac{1}{q}} \right\|_{L_p} \lesssim_{p,q,r} (|m|+1)^{\frac{1}{r}} \left\| \left( \sum_{j=-1}^\infty 2^{jsq} \left| \sum_{\mu\in\nZ} \sfk_{j,\mu}(f) 1_{I_{j,\mu}} \right|^q \right)^{\frac{1}{q}} \right\|_{L_p} . \label{eq:maximalEstimateHardDirec}
\end{align}
We make the usual modification if $q = \infty$.
\end{lem}
\begin{proof}
Step one is to show, for $x\in\nR$ and $0 < r < \min\{p,q\}$, $f\in\cB$ fixed, that
\begin{align}
    \left| \sum_{\mu\in\nZ} \sfk_{j,\mu+m}(f) 1_{I_{j,\mu}}(x) \right| \leqslant (|m|+1)^{1/r} M_r \left( \sum_{\mu\in\nZ} \sfk_{j,\mu}(f) 1_{I_{j,\mu}}(\cdot) \right)(x) . \label{eq:maximalFirstGoal}
\end{align}
We will temporarily agree that $x\in I_{j,\nu}$ for $\nu\in\nZ$ and $j\geqslant -1$ fixed. Take $I = I_{j,\nu}\cup\cdots\cup I_{j,\nu+m}$. This interval has length $2^{-j_+}(|m|+1)$ for $j\geqslant -1$. So,
\begin{align*}
    \left[ M_r \left( \sum_{\mu\in\nZ} \sfk_{j,\mu}(f) 1_{I_{j,\mu}}\right) \right]^r &\geqslant \frac{2^{j_+}}{|m|+1} \int_I \sum_{\mu\in\nZ} \sfk_{j,\mu}(f)^r 1_{I_{j,\mu}}(y) \, dy = \frac{2^{j_+}}{|m|+1} \left|\sum_{\mu=\nu}^{\nu+m} \sfk_{j,\mu}(f)^r 2^{-j_+} \right| \\
    &\geqslant \frac{1}{|m|+1} \sfk_{j,\nu+m}(f)^r ,
\end{align*}
which, for $x\in I_{j,\nu}$, is equal to
\begin{align*}
    \frac{1}{|m|+1} \left| \sum_{\mu\in\nZ} \sfk_{j,\mu+m}(f)^r 1_{I_{j,\mu}}(x) \right| .
\end{align*}
Equation \eqref{eq:maximalFirstGoal} is, then, just a rearrangement of this. We, therefore, know
\begin{align*}
    &\left\| \left( \sum_{j=-1}^\infty 2^{jsq} \left| \sum_{\mu\in\nZ} \sfk_{j,\mu+m}(f) 1_{I_{j,\mu}}(\cdot) \right|^q \right)^{1/q} \right\|_{L_p} \\
    &\hspace{+3cm}\leqslant (|m|+1)^{1/r} \left\| \left( \sum_{j=-1}^\infty M_r \left( 2^{js} \sum_{\mu\in\nZ} \sfk_{j,\mu}(f) 1_{I_{j,\mu}}(\cdot) \right)^q \right)^{1/q} \right\|_{L_p} .
\end{align*}
A quick application of \eqref{eq:BetterFeffermanStein} gives the right hand side of \eqref{eq:maximalEstimateHardDirec}. We have written the estimate above for $q < \infty$, but the analogous equation for $q = \infty$ is also valid and we can again apply \eqref{eq:BetterFeffermanStein}.
\end{proof}

We conclude this section with the Peetre maximal function. It is given, for $f\in L_1^{\text{loc}}(\nR)$, by
\begin{align}
    \Mfk^{**}_{k,A} f(x) = \esssup_{y\in\nR} \frac{|f(y)|}{(1 + 2^k |x-y|)^A} , \quad x\in\nR \label{eq:PeetreMaximal}
\end{align}
for $A > 0$ (to be chosen), $k\in\nN_0$, and $f\in \cE(t)$ for $t>0$, which we define presently: for $t > 0$,
\begin{align*}
    \cE(t) = \{f\in L_1^{\text{loc}}(\nR) : \supp \widehat{f} \subset \{|\xi| < 2t\}\} .
\end{align*}
We remark that it is more common in the literature to define $\cE(t)$ as a subset of $\cS'$ rather than $L_1^{\text{loc}}$, but this adds some technical complication to interpreting \eqref{eq:PeetreMaximal}. We will only use the Peetre maximal function for functions in $L_1^{\text{loc}}$, so we avoid these technicalities. See \cite{park2024vector} for further reading.

We now have the analog of the Fefferman-Stein inequality for the Peetre maximal function; this originally came from \cite{peetre1975spaces}, but it is discussed in many places such as \cite[Sections 1.3 and 1.4]{triebel1983function} and \cite{park2024vector} (where one may also find further results).

\begin{lem}\label{lem:PeetreMaximalIneq}
Suppose $0<p\leqslant \infty$, $A > 1/p$, and $k\in\nN_0$. Then for $f\in \cE(2^k)$,
\begin{align}
    \|\Mfk^{**}_{k,A} f\|_{L_p} \lesssim_{p,A} \|f\|_{L_p} . \label{eq:PeetreSimple}
\end{align}
Similarly, if $0<p<\infty$, $0<q\leqslant \infty$, and $A > \max\{1/p,1/q\}$, then
\begin{align}
    \|(\Mfk^{**}_{k,A} f_k)_k\|_{L_p(\ell_q)} \lesssim_{p,q,A} \|(f_k)_k\|_{L_p(\ell_q)} \label{eq:PeetreHard}
\end{align}
for any sequence $(f_k)$ with $f_k\in\cE(2^k)$ for each $k$.
\end{lem}


\subsection{Duality Considerations}\label{sec:duality}

We now discuss some points on duality, the relevance of which are twofold. First, of technical importance, this justifies that the dual pairing $(f,\psi_{j,\mu})$ makes sense for $f\in\cB$ (and $\sfk_{j,\mu}(f)<\infty$ for each $j,\mu$). Second, this will allow us to make some comments later about the optimality of the ranges \eqref{eq:sBesovRange} and \eqref{eq:sTriebelRange}. The duality we discuss will be in the framework of $(\cS(\nR),\cS'(\nR))$. The considerations below are standard, see also \cite[Remark 4.2]{schafer2021hyperbolic}.

\begin{prop}\label{prop:BLinAspq}
Let $\phi$ be an order $m\in\nN_0$ Battle-Lemari\'e wavelet or scaling function. Suppose $0<p,q\leqslant\infty$ (requiring $p<\infty$ for the $F$-spaces). Then,
\begin{enumerate}[label=(\roman*)]
    \item $\phi\in F^s_{p,q}(\nR)$ if and only if $s < m + \frac{1}{p}$ and
    \item $\phi\in B^s_{p,q}(\nR)$ if and only if either $0<q<\infty$ and $s < m+\frac{1}{p}$ or $q=\infty$ and $s\leqslant m+\frac{1}{p}$.
\end{enumerate}
\end{prop}

This proposition has been mentioned in the literature and seems to be well-known, although we have not seen a proof of it explicitly written down. Nevertheless, in \cite[Section 2.3.1, Lemma 3]{runst2011sobolev} 
the authors prove the above result for the B-splines. Their proof works again for us in proving Proposition \ref{prop:BLinAspq}; the only additional detail that is needed is that when $\phi$ is written through polynomials of order at most $m$ on each interval $(\mu,\mu+\frac12)$ for $\mu\in\frac12 \nZ$, the leading coefficients decay exponentially by Lemma \ref{lem:psiPoly}.

We now justify that $(f,\psi_{j,\mu})$ makes sense for $f\in \cB$. Of course, we need only look at $(f,\psi)$, since if $\zeta :\nR\to\nR$ is a diffeomorphism, then $g\in B^s_{p,q}$ if and only if $g\circ \zeta \in B^s_{p,q}$ and similarly for the $F$-spaces, see \cite[Section 2.10.2]{triebel1983function}.

Let us generally clarify what is meant by $(f,g)$ for $f\in\cB$ and $g\in B^{n+1}_{1,\infty}$. First, we consider two sets of two functions $\theta_0,\theta, \omega_0,\omega \in \cS$. Analogously to how we proceeded in Section \ref{subsec:spaces}, we define $\theta_k(x) = 2^k \theta(2^k x)$ for $k\geqslant 1$; we define $\omega_k$ in the same way. Further, we require that each pair satisfies \eqref{eq:betaNonzero} for some $\varepsilon > 0$, that $\omega$ has $n+2$ vanishing moments, and
\begin{align}
    \sum_{k=0}^\infty \theta_k(\xi) \omega_k(\xi) = 1 \label{eq:alphaProduct}
\end{align}
for all $\xi\in\nR$. One can choose $\omega_0,\omega$ based on $\theta_0,\theta$ so that \eqref{eq:alphaProduct} holds similarly to how we chose $\Lambda_0,\Lambda$ in Section \ref{subsec:spaces} to guarantee \eqref{eq:LPdecomp}. See also \cite[Lemma 3.6]{ullrich2016role}.

Now, let $\Theta_k = \cF^{-1}\theta_k$ and $\Omega_k = \cF^{-1} \omega_k$ for each $k\in\nN_0$ and define, for $f\in\cB$ and $g\in B^{n+1}_{1,\infty}$,
\begin{align}
    (f,g) := \sum_{k=0}^\infty (\Theta^-_k * f , \Omega_k * g) , \label{eq:dualityDef}
\end{align}
where $\Theta^-_k(x) = \Theta_k(-x)$. Above, $(\Theta^-_k * f , \Omega_k * g)$ is just integration of $\Theta^-_k * f$ against $\Omega_k * g$. From this definition,
\begin{align}
    |(f,g)| &\leqslant \sum_{k=0}^\infty |(\Theta_k^- * f , \Omega_k * g)| \leqslant \sup_{k\in\nN_0} 2^{k(n+1)} \|\Omega_k * g\|_{L_1} \sum_{k=0}^\infty 2^{k(-(n+1))} \|\Theta_k^- * f\|_{L_\infty} \notag \\
    &\lesssim \|f\|_{B^{-(n+1)}_{\infty,1}} \|g\|_{B^{n+1}_{1,\infty}} \label{eq:dualityBound} .
\end{align}
The final inequality follows simply by the definition of the Besov spaces. Since $\psi\in B^{n+1}_{1,\infty}$ by Proposition \ref{prop:BLinAspq}, the pairing $(f,\psi)$ makes sense and is finite.

In summary, it will never be incorrect in this paper to interpret the notation $(\cdot,\cdot)$ with \eqref{eq:dualityDef}. However, there are a number of situations, in which \eqref{eq:dualityDef} is equal to a much simpler expression, such as the usual pairing of a distribution with a function, e.g in Section \ref{sec:edgeCase} or Lemmata \ref{lem:scale1} and \ref{lem:scale2} (where we can even interpret $(\cdot,\cdot)$ as an integral, because we are dealing with $L_1^{\text{loc}}$ functions, one of which decays rapidly), or as an $L_2$ inner product (e.g. in the proof of Lemma \ref{lem:hardPartCrux}), since our wavelets and scaling functions are real-valued.

\begin{rem}
The definition \eqref{eq:dualityDef} actually makes sense for all $f\in\cS'$, which begs the question: why not just state our theorems for $f\in\cS'$? The answer is, we could, but it is more natural to look at $f\in\cB$ because of \eqref{eq:dualityBound}.
\end{rem}

\begin{rem}
Consider $\mathring{B}^{-(n+1)}_{\infty,1}(\nR)$, the completion of $\cS$ in $B^{-(n+1)}_{\infty,1}(\nR)$ endowed with the usual $B^{-(n+1)}_{\infty,1}(\nR)$ norm. If we were to assume $\max\{p,q\} < \infty$ in Theorems \ref{thm:mainBesov} and \ref{thm:mainTriebel}, then we could replace $\cB$ with $\mathring{B}^{-(n+1)}_{\infty,1}$ as our reference space and simply interpret $(f,\psi)$ for $f\in \mathring{B}^{-(n+1)}_{\infty,1}$ as a dual pairing. This is because, by \cite[Section 2.5.1, Remark 7]{triebel1978spaces}, we have $(\mathring{B}^{-(n+1)}_{\infty,1})' = B^{n+1}_{1,\infty}$.\footnote{We note that all of the spaces considered in Theorems \ref{thm:mainBesov} and \ref{thm:mainTriebel} with $\max\{p,q\} < \infty$ are embedded into $\mathring{B}^{-(n+1)}_{\infty,1}$ since they are embedded into $\cB$ and they only contain elements that can be approximated by Schwartz functions.}
\end{rem}

\begin{rem}
Any dependence on $\Theta_k$ and $\Omega_k$ created by defining $(f,g)$ as we did in \eqref{eq:dualityDef} is not relevant to us, because we are only evaluating these pairings in the context of a norm. Our proofs below will have no dependence on the particular choices of $\theta_k$ and $\omega_k$, so any choice of them will give an equivalent norm. We will not notate dependence on these two functions outside of this subsection.
\end{rem}



%% file: Sections/mainResultProof.tex
\section{Proofs for the Main Theorems}\label{sec:proofs}

So that we can discuss the parts of the proof below, we will say that the weaker inequality is that, which cannot be used to prove that a function is actually in $B^s_{p,q}$ or $F^s_{p,q}$ using our oversampled system, i.e. see equations \eqref{eq:easyDirecBesov} and \eqref{eq:easyDirecTriebel}. The strong inequality will, then, be defined as the one which can be used for this purpose, i.e. \eqref{eq:hardDirecBesov} and \eqref{eq:hardDirecTriebel}. These names are also justified by the proofs in the sense that in the weak direction, our proof will only depend on general properties of the Battle-Lemari\'e system (such as (I)-(IV) from Section \ref{subsec:wavelets}), whereas the strong direction will require us to look deeper at the actual structure of the system.


\subsection{The Weaker Direction}\label{sec:weak}

We start by formulating the weaker direction for each of Theorems \ref{thm:mainBesov} and \ref{thm:mainTriebel} explicitly.

\begin{prop}\label{prop:easyDirecBesov}
Suppose $p,q,s,n,\psi,\Psi$ are as in Theorem \ref{thm:mainBesov}. For each choice of $p,q,s$, there exists an $A > 1/p$ in \eqref{eq:PeetreMaximal}, such that there exists a constant $C = C(p,q,s,n,\psi,\Psi,A,\beta,\beta_0,\Lambda,\Lambda_0) > 0$ fulfilling the following:
\begin{align}
    \left\| \sfk_{j,\mu}(f) \right\|_{b^s_{p,q}} \leqslant C \|f\|_{B^s_{p,q}} \label{eq:easyDirecBesov}
\end{align}
for all $f\in\cB$.
\end{prop}

\begin{prop}\label{prop:easyDirecTriebel}
Take $p,q,s,n,\psi,\Psi$ as in Theorem \ref{thm:mainTriebel}. For each choice of $p,q,s$, there exists an $A > \max\{1/p,1/q\}$ in \eqref{eq:PeetreMaximal}, such that there exists $C = C(p,q,s,n,\psi,\Psi,A,\beta,\beta_0,\Lambda,\Lambda_0) > 0$ satisfying the following:
\begin{align}
    \left\| \sfk_{j,\mu}(f) \right\|_{f^s_{p,q}} \leqslant C \|f\|_{F^s_{p,q}} \label{eq:easyDirecTriebel}
\end{align}
for all $f\in\cB$.
\end{prop}

Below, we will make use of $\lesssim$ (or $\sim$) to avoid writing some constants, as mentioned in the Introduction. All of the parameters, upon which this $C$ in Propositions \ref{prop:easyDirecBesov} and \ref{prop:easyDirecTriebel} depends, are fixed, save $A$ and $s$, so dependency on $p,q,n,\psi,\Psi,\beta,\beta_0,\Lambda,\Lambda_0$ will not be notated.\footnote{Actually, $p,q$ are fixed differently depending on if we are dealing with $F^s_{p,q}$ or $B^s_{p,q}$, but our focus will always be clear from context.}

\begin{rem}
The dependencies of $C$ in Propositions \ref{prop:easyDirecBesov} and \ref{prop:easyDirecTriebel} on $\beta,\beta_0,\Lambda,\Lambda_0$ are superficial. They each come from either $L_p$-type bounds for these functions---e.g. in \eqref{eq:betaBound1} or the uniform bound for $F_{n+1}$ in \eqref{eq:betaBound2}---or in the conversion between the defining pieces for $B^s_{p,q}$ and $F^s_{p,q}$, see \eqref{eq:LambaBound}; this and the use of Lemma \ref{lem:PeetreMaximalIneq} are also where dependencies on $p,q,s$ come into play.
\end{rem}

\begin{rem}
We have assumed far more than we need to for the proof of Propositions \ref{prop:easyDirecBesov} and \ref{prop:easyDirecTriebel}. Indeed, we will only use properties (I)-(IV) of the Battle-Lemari\'e wavelet system in what is to follow. We will not use that the splines making up this system are orthogonal, nor that they give rise to a multiresolution analysis of $L_2(\nR)$.
\end{rem}


\subsubsection{Some Lemmata}

The main tools that we will be applying to this problem are the Littlewood-Paley pieces $(\beta_k)$ and the Peetre maximal function from Sections \ref{subsec:spaces} and \ref{subsec:maximal}, respectively. Recall our assumptions on $\beta_0,\beta$ listed directly after Remark \ref{rem:equivDef}. As always, $n\in\nN_0$ is fixed along with the particular Battle-Lemari\'e system, $\{\psi,\Psi\}$.

For this direction of the proof, oversampling is unnecessary (Propositions \ref{prop:easyDirecBesov} and \ref{prop:easyDirecTriebel} imply their weaker counterparts with no oversampling). Hence, for $p,q,s$ as in Theorem \ref{thm:mainBesov}, we want
\begin{align}
    \left( \sum_{j=-1}^\infty 2^{j\left( s-\frac{1}{p} \right)q} \left( \sum_{\mu\in\nZ} |2^j (f,\psi_{j,\mu}^\delta)|^p \right)^{q/p} \right)^{1/q} \lesssim_{s,A} \|f\|_{B^s_{p,q}} \label{eq:deltaShiftBesov}
\end{align}
or the analogous statement if $\max\{p,q\}=\infty$, and for $p,q,s$ as in Theorem \ref{thm:mainTriebel},
\begin{align}
    \left\| \left( \sum_{j=-1}^\infty 2^{jsq} \left| \sum_{\mu\in\nZ} 2^j (f,\psi_{j,\mu}^\delta) 1_{I_{j,\mu}} \right|^q \right)^{1/q} \right\|_{L_p} \lesssim_{s,A} \|f\|_{F^s_{p,q}} , \label{eq:deltaShiftTriebel}
\end{align}
making the proper edit if $q=\infty$. Above, $\psi_{j,\mu}^\delta$, for $j\in\nN_0$, $\mu\in\nZ$, is the $\delta$-shift
\begin{align*}
    \psi_{j,\mu}^\delta(x) = \psi(2^j x - (\mu+\delta))
\end{align*}
with $\delta\in [0,1)$. We are not oversampling at the base scale, so take $\psi_{-1,\mu}^\delta(x) = \Psi(x-\mu)$ for $\mu\in\nZ$. Propositions \ref{prop:easyDirecBesov} and \ref{prop:easyDirecTriebel} then follow from combining the cases $\delta = 0$ and $\delta = 1/2$.

To relate the left hand sides of \eqref{eq:deltaShiftBesov} and \eqref{eq:deltaShiftTriebel} back to the Littlewood-Paley description of the Besov and Triebel-Lizorkin scales, we first decompose $f\in\cB$ as in \eqref{eq:LPdecomp}:
\begin{align}
    f = \sum_{k=0}^\infty L_k(\Lambda_k f) = \sum_{k=0}^{\infty} L_k f_k , \label{eq:decomposef}
\end{align}
where $f_k = \Lambda_k f$ and $\widehat{f_k}$ is supported in the annulus $\{2^{k-2} \leqslant |\xi| \leqslant 2^k\}$ for $k\geqslant 1$ and in the disk $\{|\xi| \leqslant 1\}$ for $k = 0$; these support properties follow from the definition of $\Lambda_0,\Lambda$ in \eqref{eq:LambdaDef}. In particular, we will be able to apply Lemma \ref{lem:PeetreMaximalIneq} later.

We write, for $f\in \cB$ (and any $j\geqslant -1$, $\mu\in\nZ$, $\delta\in [0,1)$), using \eqref{eq:dualityDef} and \eqref{eq:decomposef},
\begin{align}
    (f,\psi_{j,\mu}^\delta) &= \sum_{\ell=0}^\infty (\Theta_\ell^- * f , \Omega_\ell * \psi_{j,\mu}^\delta) = \sum_{\ell=0}^\infty \left( \Theta_\ell^- * \left( \sum_{k=0}^\infty L_k f_k \right) , \Omega_\ell * \psi_{j,\mu}^\delta \right) \notag \\
    &= \sum_{k=0}^\infty \sum_{\ell=0}^\infty (\Theta_\ell^- * (L_k f_k) , \Omega_\ell * \psi_{j,\mu}^\delta) = \sum_{k=0}^\infty \sum_{\ell=0}^\infty (\Theta_\ell^- * f_k , \Omega_\ell * (L_k^-\psi_{j,\mu}^\delta)) \notag \\
    &= \sum_{k=0}^\infty (f_k , L_k^- \psi_{j,\mu}^\delta) . \label{eq:s<0extrabit}
\end{align}
The above steps are justified by basic properties of convolving distributions with Schwartz functions, importantly that $L_k f_k\in C^{\infty}$ and is of at most polynomial growth \cite[Section 2.2.1]{triebel1979interp}; we also use that $\psi,\Psi$ have exponential decay. By $L_k^-$ we mean convolution with $\beta_k(-\cdot)$; from here on, we will suppress this, because $\beta_k(-\cdot)$ would work just as well to define $B^s_{p,q}$ and $F^s_{p,q}$.

We have, thus, put ourselves into a situation where we want estimates for expressions of the form $|2^j(g,L_k \psi_{j,\mu}^\delta)|$. We will do this in two cases: when $j\geqslant k$ and when $j < k$.

\begin{lem}\label{lem:scale1}
Let $\{\psi,\Psi\}$ be a Battle-Lemari\'e system of order $n\in\nN_0$. Suppose $g\in L_1^{\emph{loc}}(\nR)$, $j \geqslant k$ with $j,k\in\nN_0$, $\delta\in [0,1)$, $A>0$, and $x\in I_{j,\mu}$ for some fixed $\mu\in\nZ$. Then
\begin{align*}
    |2^j (g,L_k\psi_{j,\mu}^\delta)| \lesssim_A \left( 2^{-(j-k)(n+1)} + e^{-2^{j-k}} \right) \Mfk_{k,A}^{**}g (x) .
\end{align*}
There is no dependence of the above inequality on $\mu$.
\end{lem}
\begin{proof}
Say, for now, $k\geqslant 1$ and observe that
\begin{align*}
    |2^j (g,L_k\psi_{j,\mu}^\delta)| &\leqslant 2^j \int \frac{|g(z)|}{(1 + 2^k|x-z|)^A} (1 + 2^k|x-z|)^A |\beta_k * \psi_{j,\mu}^\delta(z)| \, dz \\
    &\leqslant 2^j \Mfk_{k,A}^{**}g(x) \int (1 + 2^k|x-z|)^A |\beta_k * \psi_{j,\mu}^\delta(z)| \, dz \\
    &= 2^j \Mfk_{k,A}^{**}g(x) \int (1 + 2^k|x-z|)^A |\beta * \psi_{j-k,\mu}^\delta(2^k z)| \, dz \\
    &= 2^{j-k} \Mfk_{k,A}^{**}g(x) \int (1 + |2^k x-z|)^A |\beta * \psi_{j-k,\mu}^\delta(z)| \, dz .
\end{align*}
We analyze the above integral by splitting it in two: roughly, for now, near where $\beta * \psi_{j-k,\mu}^\delta$ achieves its maximum and otherwise. To reduce the clutter, we denote $\ell = j-k$. First, we show
\begin{align}
    \sup_{z\in\nR} |\beta * \psi_{\ell,\mu}^\delta(z)| \lesssim 2^{-\ell(n+2)} . \label{eq:globalBound1}
\end{align}
To do this, we borrow some ideas from \cite[Lemma 1]{rychkov1999theorem}. Note,
\begin{align}
    \sup_{z\in\nR} |\beta * \psi_{\ell,\mu}^\delta(z)| = \sup_{z\in\nR} |\cF^{-1}(\widehat{\beta} \widehat{\psi_{\ell,\mu}^\delta})(z)| \leqslant \|\widehat{\beta} \widehat{\psi_{\ell,\mu}^\delta}\|_{L_1} . \label{eq:firstBound}
\end{align}
Using basic properties of the Fourier transform,
\begin{align*}
    \widehat{\psi_{\ell,\mu}^\delta}(\xi) = 2^{-\ell} e^{-i (\mu + \delta)2^{-\ell}\xi} \widehat{\psi}(2^{-\ell} \xi) .
\end{align*}
Due to the $n+1$ vanishing moments of $\psi$, we can apply an argument with Taylor expansions to show
\begin{align}
    |\widehat{\psi_{\ell,\mu}^\delta}(\xi)| = |2^{-\ell} \widehat{\psi}(2^{-\ell}\xi)| \lesssim 2^{-\ell} (2^{-\ell}|\xi|)^{n+1} . \label{eq:TaylorBound}
\end{align}
The bound above is obvious for $|\xi|\geqslant 1$, because $\widehat{\psi}$ is smooth and decays to zero as $|\xi|$ approaches infinity (from basic properties of the Fourier transform and the fact that $\psi$ decays exponentially). Therefore, we must only give the bound for $|\xi| < 1$, which follows from Taylor's theorem \cite[Theorem 2.5.6]{conway2017first}. Precisely, Taylor's theorem says there exists a function, $h_{n+1} : \nR\to \nR$, such that
\begin{align}
    \widehat{\psi}(\xi) = \sum_{\kappa=0}^n \frac{\widehat{\psi}^{(\kappa)}(0)}{\kappa!} \xi^\kappa + h_{n+1}(\xi) \xi^{n+1} , \quad \xi\in (-1,1) . \label{eq:Taylor}
\end{align}
Furthermore, for each $\xi\in (-1,1)$, there exists $d$ between $\xi$ and 0 such that 
\begin{align*}
    h_{n+1}(\xi) = \frac{\widehat{\psi}^{(n+1)}(d)}{(n+1)!} .
\end{align*}
This implies that $h_{n+1}$ is uniformly bounded on $(-1,1)$, say by $C > 0$. Furthermore, due to the vanishing moments of $\psi$, we know that $\widehat{\psi}^{(\kappa)}(0) = 0$ for $0\leqslant \kappa\leqslant n$. Hence, \eqref{eq:Taylor} reveals that
\begin{align*}
    |\widehat{\psi}(\xi)| = \left| \sum_{\kappa=0}^n \frac{\widehat{\psi}^{(\kappa)}(0)}{\kappa!} \xi^\kappa + h_{n+1}(\xi) \xi^{n+1} \right| = |h_{n+1}(\xi) \xi^{n+1}| \leqslant C |\xi|^{n+1}
\end{align*}
for $\xi\in (-1,1)$. This completes the bound of $\widehat{\psi}$ near zero, so we have the bound on $\nR$ as well.

Inserting \eqref{eq:TaylorBound} in \eqref{eq:firstBound},
\begin{align}
    \sup_{z\in\nR} |\beta * \psi_{\ell,\mu}^\delta(z)|\leqslant \|\widehat{\beta} \widehat{\psi_{\ell,\mu}^\delta}\|_{L_1} \lesssim 2^{-\ell(n+2)} \|\widehat{\beta} |\xi|^{n+1}\|_{L_1} . \label{eq:betaBound1}
\end{align}
Since $\widehat{\beta}\in\cS$, the $L_1$-norm above is bounded by a constant only depending on $\beta$ and $n$. This finishes the proof of \eqref{eq:globalBound1}.

Now, we do not have compact support of our wavelet, so the term $(1 + |2^k x - z|)^A$ poses a convergence problem. To resolve this, we use the above supremum bound where $z$ is near $2^k x$ and then a bound using the exponential decay of $\psi$ away from $2^k x$ (recall $x\in I_{j,\mu}$, so this is confined to a small interval). Using that $\supp \beta \subset (-1/4,1/4)$ and property (III) of the spline,
\begin{align*}
    |\beta * \psi_{\ell,\mu}^\delta(z)| &= \left| \int_\nR \beta(z-y) \psi_{\ell,\mu}^\delta(y) \, dy \right| \leqslant \|\beta\|_\infty \int_{z-\frac14}^{z+\frac14} | \psi_{\ell,\mu}^\delta(y) | \, dy \lesssim \int_{z-\frac14}^{z+\frac14} e^{-\gamma |2^\ell y - (\mu + \delta)|} \, dy .
\end{align*}
Obviously, we then have
\begin{align}
    |\beta * \psi_{\ell,\mu}^\delta(z)| \lesssim \begin{cases}
        e^{-\gamma |2^\ell (z - \frac14) - (\mu + \delta)|} & z > 2^{-\ell}(\mu + \delta) + \frac14 , \\
        e^{-\gamma |2^\ell (z + \frac14) - (\mu + \delta)|} & z < 2^{-\ell}(\mu + \delta) - \frac14 .
    \end{cases} \label{eq:outsideBoundConv}
\end{align}
With this, we write
\begin{align*}
    \int_\nR (1 + |2^k x - z|)^A &|\beta * \psi_{\ell,\mu}^\delta(z)| \, dz = \int_{2^{-\ell}(\mu + \delta) - \frac14 - \frac{2}{\gamma}}^{2^{-\ell}(\mu + \delta) + \frac14 + \frac{2}{\gamma}} (1 + |2^k x - z|)^A |\beta * \psi_{\ell,\mu}^\delta(z)| \, dz \\
    &+ \int_{(-\infty,2^{-\ell}(\mu + \delta) - \frac14 - \frac{2}{\gamma}) \cup (2^{-\ell}(\mu + \delta) + \frac14 + \frac{2}{\gamma},\infty)} (1 + |2^k x - z|)^A |\beta * \psi_{\ell,\mu}^\delta(z)| \, dz .
\end{align*}
The expansion by $2/\gamma$ is there to let the exponential decay kick in. Within the first integral,
\begin{align*}
    |2^k x - z| \leqslant \frac{2}{\gamma} + 2 ,
\end{align*}
since $x\in I_{j,\mu}$. We also use the bound \eqref{eq:globalBound1} in this integral to see
\begin{align*}
    \left| \int_{2^{-\ell}(\mu + \delta) - \frac14 - \frac{2}{\gamma}}^{2^{-\ell}(\mu + \delta) + \frac14 + \frac{2}{\gamma}} (1 + |2^k x - z|)^A |\beta * \psi_{\ell,\mu}^\delta(z)| \, dz \right| \lesssim_A 2^{-\ell (n+2)} .
\end{align*}
Finally, the latter integral: we focus on when $z < 0$, since the part where $z > 0$ is essentially the same. Consulting \eqref{eq:outsideBoundConv} and applying the transformation $2^\ell (z + \frac14) - (\mu + \delta)\mapsto z$ reveals
\begin{align*}
    \int_{-\infty}^{2^{-\ell}(\mu + \delta) - \frac14 - \frac{2}{\gamma}} (1 + &|2^k x - z|)^A |\beta * \psi_{\ell,\mu}^\delta(z)| \, dz \leqslant 2^{-\ell} \int_{-\infty}^{-\frac{2}{\gamma}\cdot 2^\ell} (3 + |2^k x - 2^{-\ell}(z + \mu)|)^A e^{-\gamma |z|} \, dz .
\end{align*}
We have also used the triangle inequality and bounded both $2^{-\ell}\delta$ and $\frac14$ above by $1$ to consolidate terms. Within the range of integration above,
\begin{align*}
    |2^k x - 2^{-\ell}(z + \mu)| \leqslant |2^k x - 2^{-\ell}\mu| + 2^{-\ell}|z|\leqslant 1 + |z| ,
\end{align*}
since $x\in I_{j,\mu}$, so we have the bound
\begin{align*}
    2^{-\ell} \int_{-\infty}^{-\frac{2}{\gamma}\cdot 2^\ell} (4 + |z|)^A e^{-\gamma |z|} \, dz .
\end{align*}
Let us proceed crudely:
\begin{align*}
    2^{-\ell} \int_{-\infty}^{-\frac{2}{\gamma}\cdot 2^\ell} (4 + |z|)^A e^{-\gamma |z|} \, dz \leqslant 2^{-\ell} e^{-2^\ell} \int_{-\infty}^0 (4 + |z|)^A e^{-(\gamma/2) |z|} \, dz \lesssim_A 2^{-\ell} e^{-2^\ell} .
\end{align*}
The integral over the positive $z$ values may be bounded in the same way.

Turning back to the calculation we did at the very beginning of the proof, for $x\in I_{j,\mu}$ and $k\geqslant 1$,
\begin{align*}
    |2^j (g, L_k\psi_{j,\mu}^\delta)| &\lesssim_A 2^{j-k} \Mfk_{k,A}^{**} g(x) \left( 2^{-{(j-k)(n+2)}} + 2^{k-j} e^{- 2^{j-k}} \right) \\
    &= \Mfk_{k,A}^{**}g(x) \left( 2^{-(j-k)(n+1)} + e^{-2^{j-k}} \right) ,
\end{align*}
as we wished.

Notice that the above proof made no special use of the properties of $\beta$ other than its inclusion in $\cS$ and its compact support. Since $\beta_0$ also enjoys these properties, we also have the result for $k=0$.
\end{proof}

Now we do the case where $k > j$. If we were to try the same method as in Lemma \ref{lem:scale1} with the roles of $\beta$ and $\psi$ reversed, we would need to bound $\|\widehat{\psi} |\xi|^{n+1}\|_{L_1}$, which is not obviously finite, since $\psi\in C^{n-1}$ only corresponds to decay at a rate of $1 / |\xi|^{n-1}$ in phase space. We will, therefore, need a different method to bound the convolution; we do this below by building off of an idea from \cite[Lemma 6.2]{srivastava2023orthogonal}.

\begin{lem}\label{lem:scale2convolution}
Suppose $\{\psi,\Psi\}$ is a Battle-Lemari\'e system of order $n\in\nN_0$. Let $\nu\in\nZ$ and define the intervals
\begin{align}
    J_{\nu} = \left[ \frac{\nu}{2} + \frac14 , \frac{\nu}{2}  + \frac34 \right] . \label{eq:intervalJnu}
\end{align}
Then, for all $\ell\in\nN$ and $\nu\in\nZ$,
\begin{align}
    \sup_{z\in J_{\nu}} |\beta_\ell * \psi(z)| \lesssim 2^{-\ell(n+1)} e^{-(\gamma / 2)|\nu|} . \label{eq:JnuEstimate}
\end{align}
The estimate above also holds with $\Psi$ replacing $\psi$.
\end{lem}
\begin{proof}
We wish to investigate $\beta_\ell * \psi$, so we break $\psi(z-y)$ into polynomial parts, allowing us to differentiate $n$, rather than only $n-1$, times. Fix $\nu\in\nZ$ and say $z\in J_\nu$. For $\ell\in\nN$,
\begin{align}
    \supp \beta_\ell \subset \supp \beta \subset \left( z - \frac{\nu}{2} - 1 , z - \frac{\nu}{2} \right) , \label{eq:correctInterval}
\end{align}
so we can integrate over just this interval in $\beta_\ell * \psi$. Further, if $y$ is in the interval on the right-hand side of \eqref{eq:correctInterval},
\begin{align*}
    z - y \in \left( \frac{\nu}{2} , \frac{\nu}{2} + 1 \right) ,
\end{align*}
which is an interval of unit length starting on an element of $\frac12 \nZ$. Thus, using Lemma \ref{lem:psiPoly}, we can split $\psi$ into two polynomials with equal non-leading coefficients:
\begin{align*}
    \psi(z-y) = A^n_{\nu} \left( z - y - \frac{\nu + 1}{2} \right)^n + A^{n-1}\left( z - y - \frac{\nu + 1}{2} \right)^{n-1} + \cdots + A^0
\end{align*}
for $z-y \in \left[ \frac{\nu}{2} , \frac{\nu}{2} + \frac12 \right]$ and
\begin{align*}
    \psi(z-y) = A^n_{\nu + 1} \left( z - y - \frac{\nu + 1}{2} \right)^n + A^{n-1}\left( z - y - \frac{\nu + 1}{2} \right)^{n-1} + \cdots + A^0
\end{align*}
for $z-y \in \left[ \frac{\nu}{2} + \frac12 , \frac{\nu}{2} + 1 \right]$.

Splitting the interval in \eqref{eq:correctInterval} at its middle point, we may now apply our considerations to $\beta_\ell * \psi$ by replacing $\psi(z-y)$ in the integral with the above two polynomials. To save some space, we write
\begin{align*}
    z_{\nu} = z - \frac{\nu}{2} , \quad \nu\in\nZ .
\end{align*}
Now, based on what we have just discussed, for $z\in J_\nu$,
\begin{align*}
    &\beta_\ell * \psi(z) \\ &\hspace{+5mm}= A^n_{\nu + 1} \int_{z_{\nu} - 1}^{z_{\nu} - \frac12} \beta_\ell(y) \left( z - y - \frac{\nu + 1}{2} \right)^n \, dy + A^n_{\nu} \int_{z_{\nu} - \frac12}^{z_{\nu}} \beta_\ell(y) \left( z - y - \frac{\nu + 1}{2} \right)^n \, dy \\ &\hspace{+12mm}+ \sum_{i=0}^{n-1} \left( A^i \int_{z_{\nu} - 1}^{z_{\nu} - \frac12} \beta_\ell(y) \left( z - y - \frac{\nu + 1}{2} \right)^i \, dy + A^i \int_{z_{\nu} - \frac12}^{z_{\nu}} \beta_\ell(y) \left( z - y - \frac{\nu + 1}{2} \right)^i \, dy \right) .
\end{align*}
Since the non-leading coefficients are equal, also using \eqref{eq:correctInterval}, we see that for $i = 0,\ldots,n-1$
\begin{align*}
    A^i \int_{z_{\nu} - 1}^{z_{\nu} - \frac12} \beta_\ell(y) &\left( z - y - \frac{\nu + 1}{2} \right)^i \, dy + A^i \int_{z_{\nu} - \frac12}^{z_{\nu}} \beta_\ell(y) \left( z - y - \frac{\nu + 1}{2} \right)^i \, dy \\
    &= A^i \int \beta_\ell(y) \left( z - y - \frac{\nu + 1}{2} \right)^i \, dy .
\end{align*}
Recall, now, the vanishing moments of $\beta$, to see that the above integral must be zero. 

We, thus, only have the degree $n$ terms remaining: for $z\in J_\nu$,
\begin{align*}
    \beta_\ell * \psi(z) = A^n_{\nu + 1} \int_{z_{\nu} - 1}^{z_{\nu} - \frac12} \beta_\ell(y) &\left( z - y - \frac{\nu + 1}{2} \right)^n \, dy + A^n_{\nu} \int_{z_{\nu} - \frac12}^{z_{\nu}} \beta_\ell(y) \left( z - y - \frac{\nu + 1}{2} \right)^n \, dy .
\end{align*}
To see the scaling come out more clearly, we substitute $2^\ell y \mapsto y$ to get
\begin{align*}
    A^n_{\nu + 1} \int_{2^\ell (z_{\nu} - 1)}^{2^\ell(z_{\nu} - \frac12)} \beta(y) \left( z - 2^{-\ell} y - \frac{\nu + 1}{2} \right)^n \, dy + A^n_{\nu} \int_{2^\ell (z_{\nu} - \frac12)}^{2^\ell z_{\nu}} \beta(y) \left( z - 2^{-\ell} y - \frac{\nu + 1}{2} \right)^n \, dy .
\end{align*}
At this point, we integrate by parts $n$ times, for which the following notation is needed:
\begin{align*}
    F_1(y) = \int_{-\infty}^y \beta(t) \, dt , \quad F_{i+1}(y) = \int_{-\infty}^y F_i(t) \, dt 
\end{align*}
for $i\in\nN$. Due to the vanishing moments of $\beta$, each of $F_1,\ldots,F_{n+1}$ will have the same support as $\beta$, i.e. it is contained in $(-1/4,1/4)$. In any case, integration by parts yields, for $z\in J_\nu$,
\begin{align}
    \beta_\ell * \psi(z) &= A^n_{\nu + 1} 2^{-\ell n} n! \int_{2^\ell (z_{\nu} - 1)}^{2^\ell (z_{\nu} - \frac12)} 2^{-\ell} F_n(y) \, dy + A^n_{\nu} 2^{-\ell n} n! \int_{2^\ell (z_{\nu} - \frac12)}^{2^\ell z_{\nu}} 2^{-\ell} F_n(y) \, dy \notag \\
    &= 2^{-\ell (n+1)} n! (A^n_{\nu + 1} - A^n_{\nu}) F_{n+1} \left( 2^\ell \left( z - \frac{\nu}{2} - \frac12 \right) \right) , \label{eq:betaBound2}
\end{align}
where the other boundary terms have canceled out, due to the support properties of $F_{n+1}$ (and the fact that $z\in J_{\nu}$). We need two final observations: first, $F_{n+1}$ is smooth and compactly supported, so it is uniformly bounded, and second, by Lemma \ref{lem:psiPoly}, the leading coefficient terms are bounded:
\begin{align*}
    |A^n_{\nu + 1}| \leqslant 4C_0 e^{\gamma/2} e^{-(\gamma/2) |\nu + 1|} , \quad |A^n_{\nu}| \leqslant 4C_0 e^{\gamma/2} e^{-(\gamma / 2)|\nu|} .
\end{align*}
Pulling this all together, we arrive at \eqref{eq:JnuEstimate}, as we wished.

We did not make use of any property of $\psi$ that $\Psi$ does not also possess, so a simple repetition of the above steps proves the lemma for the case where $\Psi$ replaces $\psi$.
\end{proof}

Using this, the bound when $k > j$ follows easily.

\begin{lem}\label{lem:scale2}
Let $\{\psi,\Psi\}$ be a Battle-Lemari\'e system of order $n\in\nN_0$. Suppose $g\in L_1^{\emph{loc}}(\nR)$, $k > j$ with $j\geqslant -1$ and $k\in\nN_0$, $\delta\in [0,1)$, $A>0$, and $x\in I_{j,\mu}$ for $\mu\in\nZ$ fixed. Then
\begin{align*}
    |2^j (g,L_k\psi_{j,\mu}^\delta)| \lesssim_A 2^{-(k-j)(n+1)} 2^{(k-j)A} \Mfk_{k,A}^{**}g (x) .
\end{align*}
There is no dependence on $\mu$ in the above inequality.
\end{lem}
\begin{proof}
We start with a similar computation to that in Lemma \ref{lem:scale1}: for $k > j\geqslant 0$
\begin{align*}
    |2^j (g, L_k \psi_{j,\mu}^\delta)| &\leqslant 2^j \Mfk_{k,A}^{**} g(x) \int (1 + 2^k |x-z|)^A |\beta_{k-j} * \psi_{0,\mu}^\delta(2^j z)| \, dz \\
    &\leqslant \Mfk_{k,A}^{**} g(x) 2^{(k-j)A} \int (1 + |2^j x - (z+\mu+\delta)|)^A |\beta_{k-j} * \psi(z)| \, dz .
\end{align*}
We will, this time, write $\ell = k-j$. Since $x\in I_{j,\mu}$,
\begin{align*}
    |2^j x - (z+\mu+\delta)| \leqslant |2^j x - \mu| + 1 + |z| \leqslant 2 + |z| ,
\end{align*}
meaning
\begin{align*}
    |2^j (g, L_k \psi_{j,\mu}^\delta)| \lesssim \Mfk_{k,A}^{**} g(x) 2^{\ell A} \int (3 + |z|)^A |\beta_\ell * \psi(z)| \, dz .
\end{align*}
Restricting $z$ to the intervals $J_\nu$ from \eqref{eq:intervalJnu}, for $\nu\in\nZ$ , we have $3 + |z| \leqslant \frac12 (8 + |\nu|)$, so
\begin{align*}
    \int (3 + |z|)^A |\beta_\ell * \psi(z)| \, dz \lesssim_A \sum_{\nu\in\nZ} \int_{J_\nu} (8 + |\nu|)^A |\beta_\ell * \psi(z)| \, dz ,
\end{align*}
which, after using Lemma \ref{lem:scale2convolution}, leaves us with
\begin{align*}
    |2^j (g, L_k \psi_{j,\mu}^\delta)| &\lesssim_A \Mfk_{k,A}^{**} g(x) 2^{\ell A} \sum_{\nu\in\nZ} \int_{J_\nu} (8 + |\nu|)^A 2^{-\ell(n+1)}e^{-(\gamma/2)|\nu|} \, dz \\
    &\lesssim_A \Mfk_{k,A}^{**} g(x) 2^{\ell A} 2^{-\ell(n+1)} .
\end{align*}
The sum over $\nu$ will be finite and will only depend on $A$ and $\gamma$. This concludes the case $k > j \geqslant 0$.

We still need to consider what happens when $j = -1$. First, suppose $j=-1$ and $k > 0$. We have
\begin{align*}
    |2^{-1} (g,L_k\psi_{-1,\mu}^\delta)| \leqslant 2^{kA} \Mfk_{k,A}^{**}g(x) \int  (1 + |x - (z+\mu)|)^A |\beta_k * \Psi(z)| \, dz .
\end{align*}
Similarly to before, if $z\in J_\nu$ for some $\nu\in\nZ$ and $x\in I_{-1,\mu}$, then $|x - (z+\mu)| \leqslant  \frac12 \left( 4 + |\nu| \right)$, so
\begin{align*}
    |2^{-1} (g,L_k\psi_{-1,\mu}^\delta)| \leqslant 2^{kA} \Mfk_{k,A}^{**}g(x) \sum_{\nu\in\nZ} \int_{J_\nu} (6 + |\nu|)^A |\beta_k * \Psi(z)| \, dz .
\end{align*}
Applying Lemma \ref{lem:scale2convolution}, we end up at
\begin{align*}
    |2^{-1} (g,L_k\psi_{-1,\mu}^\delta)| \lesssim_A 2^{kA} \Mfk_{k,A}^{**}g(x) 2^{-k(n+1)} \lesssim_A 2^{(k+1)A} \Mfk_{k,A}^{**}g(x) 2^{-(k+1)(n+1)} ,
\end{align*}
as desired.

Finally, we consider the case when $j = -1$ and $k = 0$. Here, we cannot apply Lemma \ref{lem:scale2convolution}, but we will not need it. Explicitly, we want the bound
\begin{align*}
    |2^{-1} (g,L_0\psi_{-1,\mu}^\delta)| \lesssim 2^{-(n+1)} 2^A \Mfk_{0,A}^{**}g(x) .
\end{align*}
Following our usual steps for finding an estimate with the Peetre maximal function,
\begin{align*}
    |2^{-1}(g,L_0\psi_{-1,\mu}^\delta)| &\leqslant 2^{-1} \Mfk_{0,A}^{**} g(x) \int (1 + |x-z|)^A |(\beta_0 * \psi_{-1,\mu}^\delta) (z)| \, dz \\
    &= 2^{-1} \Mfk_{0,A}^{**}g(x) \int (1 + |x-(z+\mu)|)^A |(\beta_0 * \Psi)(z)| \, dz \\
    &\leqslant \Mfk_{0,A}^{**} g(x) \int (2 + |z|)^A |(\beta_0 * \Psi)(z)| \, dz .
\end{align*}
For the second line above, we have used the translation $z-\mu\mapsto z$ and in the final line we have used that $x\in I_{-1,\mu}$ to bound $|x-\mu|\leqslant 1$. Thus, we only need the integral to converge. Since $\beta_0$ has compact support and $\Psi$ decays exponentially, their convolution also has exponential decay, so the integral is obviously just a finite constant. We have, now, completed the proof of the lemma.
\end{proof}

When working with both Besov norms and their Triebel-Lizorkin counterparts, the work often begins with creating an estimate for the innermost terms and then creatively (more so in the Triebel-Lizorkin case) putting all of the estimates together. Let us first complete the analysis of the inner term before seeing how to combine the data.

\begin{lem}\label{lem:innerTerm}
Let $\{\psi,\Psi\}$ be a Battle-Lemari\'e system of order $n\in\nN_0$. Suppose $k\in \nN_0$, $j\in \nN_0\cup\{-1\}$, and $\ell = k-j$. We take $A>0$ and $\delta\in[0,1)$. For $f\in\cB$, let $f_k=\Lambda_k f$ if $k\geqslant 0$ and $f_k = 0$ if $k<0$. Define
\begin{align}
    a(\ell,A) = \begin{cases}
        2^{\ell(n+1)} + e^{-2^{-\ell}} & \ell \leqslant 0 , \\
        2^{(A-n-1)\ell} & \ell > 0 .
    \end{cases}\label{eq:a(lA)}
\end{align}
Then, we have the following estimate: for $x\in\nR$,
\begin{align*}
    \sum_{\mu\in\nZ} |2^j (f,\psi_{j,\mu}^\delta)| 1_{I_{j,\mu}}(x) \lesssim_A \sum_{\ell\in\nZ} a(\ell,A) \Mfk^{**}_{j+\ell,A}(f_{j+\ell})(x) .
\end{align*}
\end{lem}
\begin{proof}
Using \eqref{eq:s<0extrabit},
\begin{align*}
    \sum_{\mu\in\nZ} |2^j (f,\psi_{j,\mu}^\delta)| 1_{I_{j,\mu}}(x) \leqslant \sum_{\mu\in\nZ} \sum_{k=0}^\infty |2^j (f_k,L_k\psi_{j,\mu}^\delta)| 1_{I_{j,\mu}}(x) .
\end{align*}
For $k\geqslant 0$, we know $f_k$ is smooth by \cite[Section 2.2.1]{triebel1979interp}, so $f_k\in L_1^{\text{loc}}$ for $k\geqslant 0$ (and obviously for $k<0$). Hence, we can apply Lemmata \ref{lem:scale1} and \ref{lem:scale2} to bound the above with
\begin{align}
    \sum_{k=0}^\infty \sum_{\mu\in\nZ} 2^{-|k-j|(n+1)} 2^{A(k-j)_+} \Mfk^{**}_{k,A}(f_k)(x) 1_{I_{j,\mu}}(x) + \sum_{k\leqslant j} \sum_{\mu\in\nZ} e^{-2^{j-k}} \Mfk_{k,A}^{**}(f_k) (x) 1_{I_{j,\mu}}(x) . \label{eq:mainSum}
\end{align}
In \eqref{eq:mainSum}, the only dependency of the summands on $\mu$ is in the indicator function, so this all sums to $1_{\nR}$. We then arrive at the bound we wanted.
\end{proof}


\subsubsection{Proving the Weaker Inequalities}

In this section, we will prove Propositions \ref{prop:easyDirecBesov} and \ref{prop:easyDirecTriebel}, which is largely a matter of applying Lemma \ref{lem:innerTerm} and then Lemma \ref{lem:PeetreMaximalIneq}. The processes below are heavily inspired by \cite{garrigos2023haar}.

\begin{proof}[Proof of Proposition \ref{prop:easyDirecTriebel}]
We will be using the notation from Lemma \ref{lem:innerTerm}. The expressions in the proof below are written for the case $q < \infty$. The proof for $q=\infty$ can be given simply by making the usual modifications to the $\ell_q$ norm in each line; we omit this repetition.

Recall that it suffices to show \eqref{eq:deltaShiftTriebel}; Lemma \ref{lem:innerTerm} shows
\begin{align*}
    \left\| \left( \sum_{j=-1}^\infty 2^{jsq} \left| \sum_{\mu\in\nZ} 2^j (f,\psi_{j,\mu}^\delta) 1_{I_{j,\mu}} \right|^q \right)^{\frac{1}{q}} \right\|_{L_p} \lesssim_A \left\| \left( \sum_{j=-1}^\infty 2^{jsq} \left( \sum_{\ell\in\nZ} a(\ell,A) \Mfk^{**}_{j+\ell,A}(f_{j+\ell}) \right)^q \right)^{\frac{1}{q}} \right\|_{L_p} .
\end{align*}
Then, applying the $u$-triangle inequality with $u = \min\{p,q,1\}$, we dominate the expression above with (a constant multiple of)
\begin{align*}
    &\left(\sum_{\ell\in\nZ} \left[ a(\ell,A) \left\| \left( \sum_{j=-1}^\infty 2^{jsq} \left| \Mfk^{**}_{j+\ell,A}(f_{j+\ell}) \right|^q \right)^{1/q} \right\|_{L_p}\right]^u\right)^{1/u} \\
    &\hspace{+3cm}\leqslant \left( \sum_{\ell\in\nZ} \left[a(\ell,A) 2^{-s\ell} \left\| \left( \sum_{j\in\nZ} \left| 2^{(j+\ell)s} \Mfk_{j+\ell,A}^{**}(f_{j+\ell}) \right|^q \right)^{1/q} \right\|_{L_p} \right]^u\right)^{1/u} .
\end{align*}
We can now apply the Peetre maximal inequality \eqref{eq:PeetreHard} so long as $A > \max\{1/p,1/q\}$, which we have assumed. This means the expression above may be estimated with
\begin{align*}
    \left( \sum_{\ell\in\nZ} \left[ a(\ell,A) 2^{-s\ell} \left\| \left( \sum_{j\in\nZ} \left|2^{(j+\ell)s} f_{j+\ell} \right|^q \right)^{1/q} \right\|_{L_p} \right]^u \right)^{1/u} \\
    &\hspace{-4cm}= \left(\sum_{\ell\in\nZ}  \left[ a(\ell,A) 2^{-s\ell} \right]^u \right)^{1/u} \|(2^{ks}\Lambda_k f)\|_{L_p(\ell_q)} .
\end{align*}
Thus, \eqref{eq:LambaBound} finishes the proof so long as the sum over $\ell$ converges. Explicitly,
\begin{align}
    \sum_{\ell\in\nZ}  \left[ a(\ell,A) 2^{-s\ell} \right]^u = \sum_{\ell \leqslant 0} \left[ 2^{\ell(n+1-s)} + 2^{-s\ell}e^{-2^{-\ell}} \right]^u + \sum_{\ell > 0} \left[ 2^{\ell(A-n-1-s)} \right]^u . \label{eq:sumConvergence}
\end{align}
The former sum converges so long as $n + 1 - s > 0$ and the latter so long as $A-n-1-s < 0$, which, in total, gives the condition
\begin{align*}
    A-(n+1) < s < n+1 .
\end{align*}
Seeing as we have only required $A > \max\{1/p,1/q\}$ of $A$, for any fixed $s$ in the range \eqref{eq:sTriebelRange}, we can find an $A$ that allows for the inequality above. We have, therefore, proven Proposition \ref{prop:easyDirecTriebel}.
\end{proof}

This brings us to the case of the Besov scale, which will be reminiscent of what we have done with the $F$-spaces.

\begin{proof}[Proof of Proposition \ref{prop:easyDirecBesov}]
We will, again, use the notation from Lemma \ref{lem:innerTerm}. As in the proof of Proposition \ref{prop:easyDirecTriebel}, we will write each expression below assuming that $p,q<\infty$. Analogous proofs can be given for the cases where $q=\infty$ and $p<\infty$ by making the usual modifications. If $p=\infty$, we have to edit \eqref{eq:allbutp=infty} below to
\begin{align*}
    \sup_{\mu\in\nZ} |2^j (f,\psi_{j,\mu}^\delta)| \leqslant \left\| \sum_{\mu\in\nZ} |2^j (f,\psi_{j,\mu}^\delta)| 1_{I_{j,\mu}}(\cdot) \right\|_{L_\infty} .
\end{align*}
Otherwise, one makes the usual modifications to the expressions below to write a proof for $p=\infty$.

We only need to prove \eqref{eq:deltaShiftBesov}. Noticing that
\begin{align}
    \left( \sum_{j=-1}^{\infty} 2^{j\left( s - \frac{1}{p} \right)q} \left( \sum_{\mu\in\nZ} |2^j(f,\psi_{j,\mu}^\delta)|^p \right)^{\frac{q}{p}} \right)^{\frac{1}{q}} = \left( \sum_{j=-1}^{\infty} 2^{jsq} \left( 2^{-\frac{j}{p}} \left( \sum_{\mu\in\nZ} |2^j (f,\psi_{j,\mu}^\delta)|^p \right)^{\frac{1}{p}} \right)^q \right)^{\frac{1}{q}} \label{eq:easyBesovBegin}
\end{align}
and
\begin{align}
    2^{-j/p} \left( \sum_{\mu\in\nZ} |2^j (f,\psi_{j,\mu}^\delta)|^p \right)^{1/p} = \left\| \sum_{\mu\in\nZ} |2^j(f,\psi_{j,\mu}^\delta)| 1_{I_{j,\mu}}(\cdot) \right\|_{L_p} , \label{eq:allbutp=infty}
\end{align}
we may replicate what we have done in the Triebel-Lizorkin case. That is, using Lemma \ref{lem:innerTerm} and the $u$-triangle inequality with $u = \min\{p,1\}$,
\begin{align*}
    \left\| \sum_{\mu\in\nZ} |2^j(f,\psi_{j,\mu}^\delta)| 1_{I_{j,\mu}}(\cdot) \right\|_{L_p} \lesssim_A \left( \sum_{\ell\in\nZ} \left[ a(\ell,A) \|\Mfk_{j+\ell,A}^{**}f_{j + \ell}\|_{L_p} \right]^u \right)^{1/u} .
\end{align*}
Since $A > 1/p$, we bound $\|\Mfk_{j+\ell,A}^{**}f_{j + \ell}\|_{L_p} \lesssim_A \|f_{j+\ell}\|_{L_p}$ with Lemma \ref{lem:PeetreMaximalIneq} and \eqref{eq:easyBesovBegin} with
\begin{align*}
    \left[ \left( \sum_{j=-1}^\infty \left( 2^{jsu} \sum_{\ell\in\nZ} \left[ a(\ell,A) \|f_{j+\ell}\|_{L_p} \right]^u \right)^{q/u} \right)^{u/q} \right]^{1/u} = \left\| 2^{jsu} \sum_{\ell\in\nZ} \left[ a(\ell,A) \|f_{j+\ell}\|_{L_p} \right]^u \right\|_{\ell_{q/u}}^{1/u} .
\end{align*}
From here we can take $v = \min\{q/u , 1\}$ and apply the $v$-triangle inequality to bound the latter term above by (a constant multiple of)
\begin{align}
    \left( \sum_{\ell\in\nZ} \left\| \left( 2^{js} a(\ell,A) \|f_{j+\ell}\|_{L_p} \right)^u \right\|_{\ell_{q/u}}^v \right)^{1/(uv)} = \left( \sum_{\ell\in\nZ} a(\ell,A)^{uv} \left\| 2^{jsu} \|f_{j+\ell}\|_{L_p}^u \right\|_{\ell_{q/u}}^v \right)^{1/(uv)} . \label{eq:easyBesovAlmost}
\end{align}
The inner term can easily be rewritten as
\begin{align*}
    \left\| 2^{jsu} \|f_{j+\ell}\|_{L_p}^u \right\|_{\ell_{q/u}} = \left( 2^{-s\ell} \right)^u \left( \sum_{j=-1}^\infty 2^{(j+\ell)sq} \|f_{j+\ell}\|_{L_p}^q \right)^{u/q} = \left( 2^{-s\ell} \right)^u \left\|(2^{ks}\Lambda_k f) \right\|_{\ell_q(L_p)}^u ,
\end{align*}
where, in the last equality, we have taken $j+\ell = k$ and recalled that $f_{j+\ell} = 0$ by definition if $j+\ell < 0$. Finally, inserting this back into \eqref{eq:easyBesovAlmost}, the expression becomes
\begin{align*}
    \left( \sum_{\ell\in\nZ} a(\ell,A)^{uv} \left( 2^{-s\ell} \right)^{uv} \left\|(2^{ks}\Lambda_k f) \right\|_{\ell_q(L_p)}^{uv} \right)^{1/(uv)} \lesssim_s \|f\|_{B^s_{p,q}} \left( \sum_{\ell\in\nZ} \left[ a(\ell,A) 2^{-s\ell} \right]^{uv} \right)^{1/(uv)} ,
\end{align*}
where we have used that $\|(2^{ks} \Lambda_k f)\|_{\ell_q(L_p)}\lesssim_s \|f\|_{B^s_{p,q}}$, see \eqref{eq:LambaBound}. To finish, note that the final sum above converges so long as $A-(n+1) < s < n+1$, which is possible with the condition $A > 1/p$ for a fixed $s$ in the range \eqref{eq:sBesovRange}. We have, thus, completed the proof of Proposition \ref{prop:easyDirecBesov}.
\end{proof}


\subsection{The Stronger Direction}\label{sec:strong}

In this subsection we prove the reverse of the inequalities in Propositions \ref{prop:easyDirecBesov} and \ref{prop:easyDirecTriebel}.

\begin{prop}\label{prop:hardDirecBesov}
Suppose $p,q,s,n,\psi,\Psi$ are fixed as in Theorem \ref{thm:mainBesov}. Then there exists a constant $C = C(p,q,s,n,\psi,\Psi) > 0$, such that
\begin{align}
    \|f\|_{B^s_{p,q}} \leqslant C \left\| \sfk_{j,\mu}(f) \right\|_{b^s_{p,q}} \label{eq:hardDirecBesov}
\end{align}
for $f\in\cB$.
\end{prop}

\begin{prop}\label{prop:hardDirecTriebel}
Given $p,q,s,n,\psi,\Psi$ as in Theorem \ref{thm:mainTriebel}, there exists $C = C(p,q,s,n,\psi,\Psi) > 0$, such that
\begin{align}
    \|f\|_{F^s_{p,q}} \leqslant C \left\| \sfk_{j,\mu}(f) \right\|_{f^s_{p,q}} \label{eq:hardDirecTriebel}
\end{align}
for $f\in\cB$.
\end{prop}

In this section, constants implicitly given with the symbol $\lesssim$ (see the Introduction) may depend on $p,q,n,\psi,\Psi$ without further notation; any other dependencies, e.g. on $s$, will be clearly notated.

\subsubsection{The Relation of the B-splines to the Battle-Lemari\'e System}\label{subsubsec:remarkableConnection}

In this section, we will make the crux observation that allows us to complete this direction of the proof. This will rely on a theorem proven in \cite{ushakova2018localisation}. To state this theorem, we need to quickly summarize the efforts that lead to it.\footnote{One may check the author's master's thesis for a much more detailed summary of \cite{ushakova2018localisation}.} We fix $n\in\nN_0$ for the following discussion.

In \cite{ushakova2018localisation}, the authors examined the orthonormalization process, which took the classical B-splines (defined in Section \ref{subsec:wavelets}) and turned them into the Battle-Lemari\'e scaling function. Briefly, one may show that an orthonormal scaling function, $G$, associated to the same multiresolution analysis as the B-spline of order $n\in\nN_0$ must have the form 
\begin{align*}
    \widehat{G}(\xi) = g(\xi) \widehat{B}_n(\xi) ,
\end{align*}
where $g$ is $2\pi$-periodic and
\begin{align}
    |g(\xi)| = \left( \sum_{k\in\nZ} |\widehat{B}_n(\xi + 2\pi k)|^2 \right)^{-1/2} . \label{eq:orthonormalizationShiftFunction}
\end{align}
See \cite[Chapter 2.2]{wojtaszczyk1997mathematical} for a proof and more exposition on this fact.

The reader may have noticed that there is an obvious choice to take for $g$; this is the traditional one, as can be seen in \eqref{eq:classicalBLScaling}. However, by cleverly rewriting \eqref{eq:orthonormalizationShiftFunction}, the authors of \cite{ushakova2018localisation} managed to find $2^n$ new\footnote{When $n=0$, the process just gives back the Haar wavelet system, so in that case, there are no new wavelet systems.} Battle-Lemari\'e scaling functions and, therefore, $2^n$ new Battle-Lemari\'e wavelets as well. These were written as $\Psi_{t_1,\ldots, t_n}^\pm$ and $\psi_{t_1,\ldots, t_n}^\pm$, respectively. Each of $t_1,\ldots, t_n$ represents a binary choice between one of two real numbers: $r_j$ or $1/r_j$ for each $j=1,\ldots,n$. These are simply constants that appear in the process of orthonormalization.

In \cite{ushakova2018localisation}, it is proven that, for $t_1,\ldots,t_n$ fixed, there exists an exponentially decaying sequence, say $(d_k)\subset\nR$, such that
\begin{align*}
    \Psi_{t_1,\ldots, t_n}^\pm (x) = \sum_{k\in\nZ} d_k B_n(x-k) .
\end{align*}
Of course, the B-splines have have compact support, so the support of $\Psi_{t_1,\ldots,t_n}^\pm$ is inherited from $(d_k)$, which means $\Psi_{t_1,\ldots,t_n}^\pm$ decays exponentially. Likewise, the exponential decay of $\psi_{t_1,\ldots, t_n}^\pm$ (with $t_1,\ldots,t_n$ fixed) is seen through the two-scale relation: there exists an exponentially decaying sequence $(e_k)\subset\nR$, such that
\begin{align*}
    \psi_{t_1,\ldots, t_n}^\pm (x) = \sum_{k\in\nZ} e_k B_n(2x-k) ,
\end{align*}
so $\psi_{t_1,\ldots, t_n}^\pm$ also decays exponentially.

Due to the derivative identity \eqref{eq:magicIdentity} and a beautiful observation (see \cite[Equation (34)]{ushakova2018localisation} and the surrounding discussion), the authors of \cite{ushakova2018localisation} were able to prove that $B_{2n+1}^{(n+1)}(2x+n)$ shares localization properties with these new Battle-Lemari\'e wavelets. This is made precise by the following theorem.

\begin{lem}[Theorem 3.1 in \cite{ushakova2018localisation}]\label{lem:ushakova}
Let $n\in\nN$. The B-spline $B_{2n+1}^{(n+1)}(2x+n)$ may be written as a finite linear combination of $\psi_{t_1,\ldots,t_n}^\pm$ and their half-integer shifts, i.e. using elements coming from
\begin{align*}
    \left\{ \psi_{t_1,\ldots,t_n}^\pm \left(\cdot-\frac{\mu}{2} \right) : \mu\in\nZ , t_j = r_j \text{ or } t_j = 1/r_j \right\} .
\end{align*}
\end{lem}

\begin{rem}
In \cite[Theorem 3.1]{ushakova2018localisation}, the authors gave an explicit finite linear combination, which achieves the above result. Due to it being rather cumbersome to write down, we do not include it. With this in hand, we can make the observation that we have been working towards.
\end{rem}

\begin{lem}\label{lem:hardPartCrux}
Let $n\in\nN_0$ and say $\psi$ is a Battle-Lemari\'e wavelet of order $n$. We can write
\begin{align}
    B_{2n+1}^{(n+1)}(2x) = \sum_{\ell\in\nZ} q_\ell \widetilde{\psi}_{0,\ell}(x) , \label{eq:hardCruxEquation}
\end{align}
where $(q_\ell)$ decays exponentially.
\end{lem}
\begin{proof}
The case $n=0$ is not included in Lemma \ref{lem:ushakova}, but this is because it is obvious and well-known. When $n=0$, due to the derivative identity \eqref{eq:magicIdentity},
\begin{align}
    B_1'(2x) = B_0(2x) - B_0(2x-1) = 1_{[0,1/2)}(x) - 1_{[1/2,1)}(x) = h(x) . \label{eq:Haar2scale}
\end{align}
The function, $h$, above is the Haar function, i.e. the classical Battle-Lemari\'e wavelet of order zero, see \eqref{eq:HaarSystem}. Hence, in this case, we can take $q_\ell = 1$ for $\ell=0$ and $q_\ell=0$ otherwise in \eqref{eq:hardCruxEquation}, which clearly decays exponentially. Recalling Remark \ref{rem:Haar}, any other Battle-Lemari\'e wavelet of order zero, call it $\psi$, satisfies $\psi(x) = h(x-k)$ for some $k\in\nZ$. In other words, for this $\psi$, we can take $q_\ell = 1$ for $\ell = 2k$ and $q_\ell = 0$ otherwise. This establishes the case $n=0$.

We, therefore, focus on when $n\in\nN$, where we can apply Lemma \ref{lem:ushakova}. As a small convenience, we prove the lemma for $B_{2n+1}^{(n+1)}(2x+n)$ instead of $B_{2n+1}^{(n+1)}(2x)$. Note that we could split the finite linear combination giving $B_{2n+1}^{(n+1)}(2x+n)$ from Lemma \ref{lem:ushakova} into two parts: that made up of functions in
\begin{align*}
    \left\{ \psi_{t_1,\ldots,t_n}^\pm(x-\mu) : \mu\in\nZ , t_j = r_j \text{ or } t_j = 1/r_j \right\} ,
\end{align*}
call this part $g_0$, and the other part with functions in
\begin{align*}
    \left\{ \psi_{t_1,\ldots,t_n}^\pm\left(x-\frac{2\mu+1}{2} \right) : \mu\in\nZ , t_j = r_j \text{ or } t_j = 1/r_j \right\} ,
\end{align*}
call this part $g_1$. So $B_{2n+1}^{(n+1)}(2x+n) = g_0(x) + g_1(x)$. Observe that the function $g_0$ is, by definition, a finite linear combination of other functions, all of which have exponential decay. Therefore, $g_0$ must also have exponential decay. The same can be said for $g_1$. 

Let $V_0$ be the base scale space coming from $B_n$ and $W_0$ its associated wavelet space; then, clearly, $g_0\in W_0$, since, by definition, Battle-Lemari\'e systems generate the same multiresolution analysis as the B-spline of the same order. We also have that $(\widetilde{\psi}_{0,2\mu})$ generates $W_0$, since $\widetilde{\psi}_{0,2\mu} = \psi_{0,\mu}$ for each $\mu\in\nZ$. Hence, since $g_0\in W_0$, there exists $(q_\ell')\in \ell_2$, such that
\begin{align}
    g_0 = \sum_{\ell\in\nZ} q_\ell' \widetilde{\psi}_{0,2\ell} . \label{eq:g0}
\end{align}
The function $g_1$ can be analyzed similarly by making a slight adjustment to our argument. Define
\begin{align*}
    W_0^{1/2} = \{f(\cdot-1/2) : f\in W_0\} .
\end{align*}
By the definition of $g_1$, it is clear that $g_1\in W_0^{1/2}$ and $(\widetilde{\psi}_{0,2\mu + 1})$ generates $W_0^{1/2}$ just as $(\widetilde{\psi}_{0,2\mu})$ generates $W_0$. Hence, there exists $(q_\ell'')\in\ell_2$ such that
\begin{align}
    g_1 = \sum_{\ell\in\nZ} q_\ell'' \widetilde{\psi}_{0,2\ell+1} . \label{eq:g1}
\end{align}
Putting these two observations about $g_0$ and $g_1$ together,
\begin{align*}
    B_{2n+1}^{(n+1)}(2x+n) = g_0(x) + g_1(x) = \sum_{\ell\in\nZ} q_\ell' \widetilde{\psi}_{0,2\ell} + \sum_{\ell\in\nZ} q_\ell'' \widetilde{\psi}_{0,2\ell+1} = \sum_{\ell\in\nZ} q_\ell \widetilde{\psi}_{0,\ell}(x) ,
\end{align*}
where $q_\ell$ is defined to be $q_\ell'$ for $\ell$ even and $q_\ell''$ for $\ell$ odd. Clearly $(q_\ell)\in \ell_2$, so it remains to demonstrate that $(q_\ell)$ decays exponentially.

Recall that $\widetilde{\psi}_{0,2\ell} = \psi_{0,\ell}$ for $\ell\in\nZ$, which is an orthonormal system. Therefore, we may take $L_2$ inner products on each side of \eqref{eq:g0} with some $\widetilde{\psi}_{0,2\ell_0}$ for $\ell_0\in\nZ$ to get
\begin{align*}
    (g_0,\widetilde{\psi}_{0,2\ell_0}) = q_{\ell_0}' .
\end{align*}
We can now estimate the inner product on the left hand side above. Since $g_0$ decays exponentially, we can write $|g_0(x)| \leqslant C_1 e^{-\gamma_1 |x|}$ for some constants $C_1,\gamma_1 > 0$. We also have, by property (III) of $\psi$, that $|\widetilde{\psi}_{0,2\ell_0}(x)| \leqslant C_0 e^{-\gamma |x-\ell_0|}$. Therefore, writing $\gamma_2 = \min\{\gamma,\gamma_1\}$,
\begin{align*}
    |q_{\ell_0}'| = |(g_0,\widetilde{\psi}_{0,2\ell_0})| \leqslant C_0 C_1 \int_{\nR} e^{-\gamma_2 |x|} e^{-\gamma_2 |x-\ell_0|} \, dx \leqslant C_0 C_1 C_{\gamma_3} e^{-\gamma_3 |\ell_0|} ,
\end{align*}
where $0<\gamma_3 < \gamma_2$ may be chosen at will and $C_{\gamma_3} > 0$ depends on $\gamma_3$. This easy calculation, which we have left out, can be found, e.g., in \cite[Lemma 3.16]{wojtaszczyk1997mathematical}.

The same process as above can now be repeated for $g_1$ to show that the sequence $(q_\ell'')_{\ell\in\nZ}\in\ell_2$ also decays exponentially. Since both $(q_\ell')$ and $(q_\ell'')$ decay exponentially, $(q_\ell)$ does too.
\end{proof}


\begin{rem}
In creating a spline frame, as we have, the language that is often used is that we are ``oversampling.'' This makes one ask the naive question of whether or not more sampling would help. Could we, e.g., sample at the integers and $1/3$-shifts? The lemma above is precisely the point in the proof where this oversampling comes into play. We, therefore, see that the need for oversampling emerges only because our proof exploits the interaction of the Battle-Lemari\'e wavelet with the B-splines, which are scaling functions. 

One could, therefore, see it as inevitable that we needed to oversample, because, if we are to use the interactions of wavelets with scaling functions, then this factor of two will always be there. This is the so-called two-scale relation and is a central idea in wavelet construction. From this perspective, we remark that with our method of proof, oversampling more, or, indeed, anywhere except at the half-integers, will not help. Additionally, since the oversampling emerges from our detour through the B-splines, this explains why we do not oversample at the base scale.
\end{rem}


\subsubsection{Setup for the Proof}\label{sec:setupStrongProof}

In this subsection, we will build the foundations for the proofs of Propositions \ref{prop:hardDirecBesov} and \ref{prop:hardDirecTriebel}. As such, we will use the assumptions (e.g. on $p,q,s,n,\psi,\Psi$) and notation involved in these two propositions throughout the subsection. 

Observe that, because of Theorem \ref{thm:SplineCharac}, Theorems \ref{thm:mainBesov} and \ref{thm:mainTriebel} are already proven for most of the range for $s$. Therefore, we look for a way to sneak through these lower parameter ranges; for that, we need \cite[Theorem 2.3.8]{triebel1983function}, which gives us, for $f\in\cB$,
\begin{align}
    \|f\|_{B^s_{p,q}(\nR)} \lesssim_s \|f\|_{B^{s-(n+1)}_{p,q}(\nR)} + \|f^{(n+1)}\|_{B^{s-(n+1)}_{p,q}(\nR)} \label{eq:BesovSumBound}
\end{align}
and
\begin{align}
    \|f\|_{F^s_{p,q}(\nR)} \lesssim_s \|f\|_{F^{s-(n+1)}_{p,q}(\nR)} + \|f^{(n+1)}\|_{F^{s-(n+1)}_{p,q}(\nR)} . \label{eq:TriebelSumBound}
\end{align}
Of course, \eqref{eq:BesovSumBound} and \eqref{eq:TriebelSumBound} are valid under the assumptions on $p,q,s,n$ in Propositions \ref{prop:hardDirecBesov} and \ref{prop:hardDirecTriebel}, respectively, but, in fact, they are valid more generally for all $0<p,q\leqslant\infty$ ($p<\infty$ for the $F$-spaces) and $s\in\nR$.

Before proving the bounds that we want, we will need to examine how $f^{(n+1)}$ behaves in the Besov and Triebel-Lizorkin norms. For a fixed $n\in\nN_0$, let us write $\phi$ to mean one of the Battle-Lemari\'e wavelets found in \cite{ushakova2018localisation} of order $2n+1$. This can be written with its two-scale identity; as mentioned in Section \ref{subsubsec:remarkableConnection}, there exists an exponentially decaying $(b_k)\subset\nR$ with
\begin{align*}
    \phi(x) = \sum_{k\in\nZ} b_k B_{2n+1} (2x - k) ,
\end{align*}
so for $j\geqslant 0$ and $\mu\in\nZ$
\begin{align}
    \phi_{j,\mu}(x) = \sum_{k\in\nZ} b_k B_{2n+1} (2^{j+1} x - (2\mu + k)) . \label{eq:2scalephi}
\end{align}
The corresponding scaling function of order $2n+1$, which we will call $\Phi$, can be written
\begin{align}
    \Phi(x) = \sum_{k\in\nZ} c_k B_{2n+1}(x-k) , \label{eq:2scalePhi}
\end{align}
where, again, $(c_k)$ decays exponentially, see Section \ref{subsubsec:remarkableConnection}. We will, as usual, write $\phi_{-1,\mu}(x) = \sqrt{2} \Phi(x-\mu)$ for $\mu\in\nZ$.

Due to Theorem \ref{thm:SplineCharac}, under the assumptions of Proposition \ref{prop:hardDirecTriebel}, we have\footnote{We can make sense of the pairing $(f^{(n+1)},\phi_{j,\mu})$ in the same way as we did for $(f,\psi_{j,\mu})$; see Section \ref{sec:duality}. Just note that, since $f\in\cB$, we know $f^{(n+1)}\in B^{-(2n+2)}_{\infty,1}$ by \cite[Theorem 2.3.8]{triebel1983function} and $\phi\in B^{2n+2}_{1,\infty}$ by Proposition \ref{prop:BLinAspq}.}
\begin{align}
    \|f^{(n+1)}\|_{F^{s-(n+1)}_{p,q}} \sim_s \left\| \left( \sum_{j=-1}^\infty 2^{j(s-(n+1))q} \left| \sum_{\mu\in\nZ} 2^j (f^{(n+1)} , \phi_{j,\mu}) 1_{I_{j,\mu}}(\cdot) \right|^q \right)^{1/q} \right\|_{L_p} , \label{eq:ChuiWangTriebel}
\end{align}
making the usual modification if $q = \infty$. To be explicit, we have assumed \eqref{eq:sTriebelRange} of $s$, which means
\begin{align}
    \max\left\{ -\frac{1}{p'} , -\frac{1}{q'} \right\} - (2n+1) < s - (n+1) < 0 . \label{eq:phisrange0}
\end{align}
By Theorem \ref{thm:SplineCharac}, the equivalence \eqref{eq:ChuiWangTriebel} is valid for
\begin{align}
    \max\left\{ -\frac{1}{p'},-\frac{1}{q'} , 0 \right\} - (2n+1) < s - (n+1) < 2n+1 \label{eq:phisrange1}
\end{align}
if $p$ or $q$ is in $(0,1]$ or $q = \infty$ and
\begin{align}
    \max\left\{ -\frac{1}{p'},-\frac{1}{q'} \right\} - (2n+1) < s - (n+1) < 2n+1 + \min\left\{ \frac{1}{p},\frac{1}{q} \right\} \label{eq:phisrange2}
\end{align}
if $1<p,q<\infty$, since $\phi$ is of order $2n+1$. Clearly \eqref{eq:phisrange0} fits comfortably both in \eqref{eq:phisrange1} and \eqref{eq:phisrange2} in each respective case.

There is a corresponding version for the Besov spaces: under the assumptions of Proposition \ref{prop:hardDirecBesov},
\begin{align}
    \|f^{(n+1)}\|_{B^{s-(n+1)}_{p,q}} \sim_s \left( \sum_{j=-1}^\infty 2^{j\left( s - (n+1) - \frac{1}{p} \right)q} \left( \sum_{\mu\in\nZ} |2^j(f^{(n+1)},\phi_{j,\mu})|^p \right)^{q/p} \right)^{1/q} , \label{eq:ChuiWangBesov}
\end{align}
making the usual modifications if $\max\{p,q\} = \infty$. Again, to be explicit, we have assumed \eqref{eq:sBesovRange} of $s$, which means
\begin{align}
    -\frac{1}{p'} - (2n+1) < s-(n+1) < 0 . \label{eq:phisrangeB0}
\end{align}
By Theorem \ref{thm:SplineCharac}, the equivalence \eqref{eq:ChuiWangBesov} holds for
\begin{align}
    -\frac{1}{p'} - (2n+1) < s-(n+1) < (2n+1) + \min\left\{ \frac{1}{p} , 1\right\} , \label{eq:phisrangeB1}
\end{align}
because $\phi$ is of order $2n+1$. Since \eqref{eq:phisrangeB0} sits in \eqref{eq:phisrangeB1}, the equivalence \eqref{eq:ChuiWangBesov} is valid.

We, thus, once again, wish to bound the innermost part, which, this time, is $2^j |(f^{(n+1)},\phi_{j,\mu})|$. We will show the following lemma, which is the crux observation allowing us to prove Propositions \ref{prop:hardDirecBesov} and \ref{prop:hardDirecTriebel}.

\begin{lem}\label{lem:hardDirecInnerBound}
Take $n\in\nN_0$ and, as above, let $\{\phi,\Phi\}$ be one of the Battle-Lemari\'e systems of order $2n+1$ from \cite{ushakova2018localisation}. As usual, $\{\psi,\Psi\}$ is a Battle-Lemari\'e system of order $n$. Below, the shifted spline coefficients are defined through $\psi,\Psi$ (not $\phi,\Phi$) as in \eqref{eq:shiftedSplineCoef} and \eqref{eq:shiftedSplinej=-1}. Suppose $f\in\cB$. With $j\geqslant -1$, $\mu\in\nZ$, we have
\begin{align}
    2^j |(f^{(n+1)}, \phi_{j,\mu})| \leqslant 2^{(j+1)(n+1)} \sum_{m\in\nZ} r_m \sfk_{j,\mu+m}(f) ,\label{eq:hardExactlyWhatWeNeed}
\end{align}
for an exponentially decaying sequence, $(r_m)$.
\end{lem}
\begin{proof}
For this proof, we will suppress the use of \eqref{eq:dualityDef} as its inclusion will not change the process below---it will only make the steps appear more complicated. Operations such as moving derivatives between components in $(\cdot,\cdot)$ continue to be permissible, as one can easily show.

For now, we let $j\geqslant 0$ and fix $\mu\in\nZ$ as well as $f\in\cB$; using the two-scale identity relating $\phi_{j,\mu}$ to $B_{2n+1}$, i.e. \eqref{eq:2scalephi}, yields
\begin{align}
    |(f^{(n+1)} , \phi_{j,\mu})| &\leqslant \sum_{k\in\nZ} |b_k| \left|(f^{(n+1)} , B_{2n+1}(2^{j+1}\cdot - (2\mu + k))) \right| \notag \\
    &= \left| (-1)^{n+1} 2^{(j+1)(n+1)} \right| \sum_{k\in\nZ} |b_k| \left|(f, B_{2n+1}^{(n+1)}(2^{j+1}\cdot - (2\mu + k))) \right| \label{eq:step1} .
\end{align}
Using Lemma \ref{lem:hardPartCrux} to rewrite the B-spline, for $\nu\in\nZ$,
\begin{align*}
    B_{2n+1}^{(n+1)}(2^{j+1}x - \nu) = \sum_{\ell\in\nZ} q_\ell \widetilde{\psi}_{j,\ell+\nu}
\end{align*}
with $(q_\ell)$ decaying exponentially, meaning that
\begin{align}
    \left| (f, B_{2n+1}^{(n+1)}(2^{j+1}\cdot - \nu)) \right| \leqslant \sum_{\ell\in\nZ} |q_\ell| \left|(f , \widetilde{\psi}_{j,\ell+\nu}) \right| . \label{eq:step2}
\end{align}
Combining \eqref{eq:step1} with \eqref{eq:step2}, we see that for $j\geqslant 0$ and $\mu\in\nZ$,
\begin{align*}
    2^j |(f^{(n+1)}, \phi_{j,\mu})| &\leqslant 2^{(j+1)(n+1)} \sum_{k\in\nZ} |b_k| 2^j \left|(f,B_{2n+1}^{(n+1)}(2^{j+1}\cdot - (2\mu + k)))\right| \\
    &\leqslant 2^{(j+1)(n+1)} \sum_{k\in\nZ} |b_k| \sum_{\ell\in\nZ} |q_\ell| 2^j \left|(f,\widetilde{\psi}_{j,\ell + (2\mu + k)})\right| .
\end{align*}
Apply the translation $m = k+\ell$:
\begin{align}
    2^j |(f^{(n+1)},\phi_{j,\mu})| \leqslant 2^{(j+1)(n+1)} \sum_{m\in\nZ} \left[ \sum_{k\in\nZ} |b_k| |q_{m-k}| \right] 2^j |(f,\widetilde{\psi}_{j,2\mu + m})| . \label{eq:step3}
\end{align}
The sum in brackets defines a new sequence in $m$, call it $(\widetilde{q}_m)_{m\in\nZ}$. Since both $(q_m)$ and $(b_k)$ decay exponentially, this discrete convolution will as well, see \cite[Lemma 3.16]{wojtaszczyk1997mathematical}. 

Recall that we are looking for a bound in terms of $\sfk_{j,\mu}(f)$; with this in mind, we can group consecutive terms by increasing some of the coefficients. That is, define
\begin{align}
    p_m = \widetilde{q}_{2m} + \widetilde{q}_{2m+1} , \label{eq:pmdef}
\end{align}
so that we may write, using \eqref{eq:step3},
\begin{align*}
    2^j |(f^{(n+1)},\phi_{j,\mu})| &\leqslant 2^{(j+1)(n+1)} \sum_{m\in\nZ} p_m \left(2^j |(f,\widetilde{\psi}_{j,2\mu + 2m})| + 2^j |(f,\widetilde{\psi}_{j,2\mu + 2m+1})| \right) \\
    &= 2^{(j+1)(n+1)} \sum_{m\in\nZ} p_m \sfk_{j,\mu + m}(f)
\end{align*}
for $\mu\in\nZ$ and $j\geqslant 0$. Importantly, the sequence $(p_m)$ does not depend in any way on $\mu$.

Now we return to the case $j=-1$. This is, in fact, mostly a repetition of the above process in an easier setting. Note that, for a fixed $\mu\in\nZ$, using \eqref{eq:2scalePhi},
\begin{align}
    2^{-1} |(f^{(n+1)},\phi_{-1,\mu})| \leqslant \frac{\sqrt{2}}{2} \sum_{k\in\nZ} |c_k| \left|(f,B_{2n+1}^{(n+1)}(\cdot - (\mu+k)))\right| . \label{eq:step4}
\end{align}
Due to the derivative identity \eqref{eq:magicIdentity}, we may collapse $B_{2n+1}^{(n+1)}(x)$ into a finite linear combination of $B_n(x)$ and its integer translates. It therefore lies in $V_0$, the scale space corresponding to $B_n$, which is generated by $(\widetilde{\psi}_{-1,m})_{m\in\nZ}$. Consequently, we can find a sequence, $(\widetilde{p}_m)\in \ell_2$, with
\begin{align}
    B_{2n+1}^{(n+1)}(x) = \sum_{m\in\nZ} \widetilde{p}_m \widetilde{\psi}_{-1,m}(x) . \label{eq:step5}
\end{align}
Due to the orthogonality in $L_2(\nR)$ of the system $(\widetilde{\psi}_{-1,m})$, we may easily show that  $(\widetilde{p}_m)$ decays exponentially: for each $m\in\nZ$,
\begin{align}
    \widetilde{p}_m = (B_{2n+1}^{(n+1)} , \widetilde{\psi}_{-1,m}) . \label{eq:expDecaySMT}
\end{align}
We know that $\Psi$ decays exponentially and $B_{2n+1}^{(n+1)}$ has compact support, so using \cite[Lemma 3.16]{wojtaszczyk1997mathematical}, we see from \eqref{eq:expDecaySMT} that $\widetilde{p}_m$ decays exponentially. We explained this proof of exponential decay more fully in the proof of Lemma \ref{lem:hardPartCrux}.

Hence, combining \eqref{eq:step4} with \eqref{eq:step5}, we see, for $\mu\in\nZ$,
\begin{align*}
    2^{-1} |(f^{(n+1)},\phi_{-1,\mu})| \leqslant \frac{\sqrt{2}}{2} \sum_{k\in\nZ} |c_k| \sum_{m\in\nZ} |\widetilde{p}_{m-k}| \left|(f, \widetilde{\psi}_{-1,\mu+m})\right| .
\end{align*}
Once again, since $(c_k)$ and $(\widetilde{p}_m)$ decay exponentially, we have that the new sequence given by
\begin{align*}
    \widetilde{r}_m = \sum_{k\in\nZ} |c_k| |\widetilde{p}_{m-k}|
\end{align*}
decays exponentially by \cite[Lemma 3.16]{wojtaszczyk1997mathematical}. Hence, for $\mu\in\nZ$,
\begin{align*}
    2^{-1} |(f^{(n+1)},\phi_{-1,\mu})| \leqslant \sum_{m\in\nZ} \widetilde{r}_m \sfk_{-1,\mu+m}(f) ,
\end{align*}
which is exactly what we wanted. Since both $(\widetilde{r}_m)$ and $(p_m)$ (from \eqref{eq:pmdef}) decay exponentially, we can just create a new sequence defined by $r_m = \max\{\widetilde{r}_m,p_m\}$ for each $m\in\nZ$. Obviously, this new sequence also decays exponentially and fulfills \eqref{eq:hardExactlyWhatWeNeed}, completing the proof.
\end{proof}

\subsubsection{Proving the Stronger Direction}\label{subsubsec:proveHardDirec}

The time has come for our work to culminate in a proof for Propositions \ref{prop:hardDirecBesov} and \ref{prop:hardDirecTriebel}, so we will continue to tacitly assume the requirements of those propositions in this subsection. Recall that we wanted to bound each of the terms in \eqref{eq:BesovSumBound} and \eqref{eq:TriebelSumBound}. The process we follow for the first term in these equations is only an adaptation of that undertaken in \cite[Proposition 4.9]{garrigos2023haar}.

\begin{lem}\label{lem:hardDirecFirstTerm}
Fix $n\in\nN_0$ and let $\{\psi,\Psi\}$ be a Battle-Lemari\'e wavelet system of order $n$. Take $0 < p,q \leqslant \infty$ and suppose
\begin{align*}
	-\frac{1}{p'} - n < s < n + 1 + \frac{1}{p} .
\end{align*}
For $f\in\cB$, we have
\begin{align*}
	\|f\|_{B^{s-1}_{p,q}} \lesssim_s \left\| (2^j (f,\psi_{j,\mu}) ) \right\|_{b^s_{p,q}} .
\end{align*}
If $p < \infty$, then for $f\in\cB$,
\begin{align*}
	\|f\|_{F^{s-1}_{p,q}} \lesssim_s \left\| (2^j (f,\psi_{j,\mu}) ) \right\|_{f^s_{p,q}} .
\end{align*}
\end{lem}
\begin{proof}
Let us focus on the Triebel-Lizorkin spaces. The range for $s$ above is more reminiscent of the normal range we have seen for the Besov spaces; that is because we will sneak through the Besov spaces. Note first that $B^r_{p,p} \subset F^{s-1}_{p,q}$ by \cite[Section 2.3.2]{triebel1983function} so long as $s-1 < r$. We examine two cases: $s > n + \frac{1}{p}$ and $s\leqslant n + \frac{1}{p}$.

In the former, we choose $r = n + \frac{1}{p} - \varepsilon$ for an $\varepsilon > 0$ chosen so that
\begin{align}
	-\frac{1}{p'} - n < r < n + \frac{1}{p} , \quad s-1 < r < s . \label{eq:howChooseEps}
\end{align}
It is possible to choose such an $\varepsilon$, because $s < n + 1 + \frac{1}{p}$, per the assumptions of this lemma. If, on the other hand, $s \leqslant n + \frac{1}{p}$, then we choose $r = s-\varepsilon$ and arrange once again for $\varepsilon > 0$ to be small enough so that \eqref{eq:howChooseEps} holds, this being possible now because $s > -\frac{1}{p'} - n$.

The reason we have required $r < s$ is so that we also have the embedding $f^s_{p,q} \subset b^r_{p,p}$ by \cite[Section 2.3.2]{triebel1983function}.\footnote{This reference gives the result for the function spaces: $F^s_{p,q}\subset B^r_{p,p}$. By Theorem \ref{thm:SplineCharac}, these function spaces are isomorphic to their sequence space counterparts, so the embedding on the side of the sequence spaces follows easily.} With this work, we have the following string of inequalities:
\begin{align*}
	\|f\|_{F^{s-1}_{p,q}} \lesssim \|f\|_{B^r_{p,p}} \lesssim \|(2^j (f,\psi_{j,\mu}))\|_{b^r_{p,p}} \lesssim \| (2^j (f,\psi_{j,\mu}) )\|_{f^s_{p,q}} .
\end{align*}
The second inequality above is just due to Theorem \ref{thm:SplineCharac}. This completes the Triebel-Lizorkin case.

We may now repeat the steps above to do the Besov case. We do not sneak through the Triebel-Lizorkin spaces, because the results for those are worse. That is, we observe with \cite[Section 2.3.2]{triebel1983function} that $B^r_{p,p} \subset B^{s-1}_{p,q}$ for $s-1 < r$ and $b^s_{p,q} \subset b^r_{p,p}$ for $r < s$ and choose $r$ just as we did above. The argument is entirely analogous.
\end{proof}

We are finally in a position to complete the proofs for Propositions \ref{prop:hardDirecBesov} and \ref{prop:hardDirecTriebel}.

\begin{proof}[Proof of Proposition \ref{prop:hardDirecTriebel}]
Recall from \eqref{eq:TriebelSumBound} that we must bound $\|f\|_{F^{s-(n+1)}_{p,q}}$ and $\|f^{(n+1)}\|_{F^{s-(n+1)}_{p,q}}$. The former can be bounded trivially using Lemma \ref{lem:hardDirecFirstTerm}. Indeed, from \cite[Section 2.3.2]{triebel1983function},
\begin{align*}
    \|f\|_{F^r_{p,q}} \lesssim_{r,\varepsilon} \|f\|_{F^{r+\varepsilon}_{p,q}}
\end{align*}
for any $r\in\nR$ and $\varepsilon > 0$, so we can take $r = s-(n+1)$ and $\varepsilon = n$. Combining this with Lemma \ref{lem:hardDirecFirstTerm} and the simple observation that
\begin{align*}
    \left\| (2^j (f,\psi_{j,\mu}) ) \right\|_{f^s_{p,q}} \leqslant \left\| \sfk_{j,\mu}(f) \right\|_{f^s_{p,q}}
\end{align*}
gives the desired estimate for the first term.

Now we look at the latter term. The steps below are written for $q < \infty$, but each can be modified in the usual way to produce a proof for the case $q=\infty$. Using \eqref{eq:ChuiWangTriebel} and Lemma \ref{lem:hardDirecInnerBound},
\begin{align*}
    \|f^{(n+1)}\|_{F^{s-{(n+1)}}_{p,q}} &\lesssim_s \left\| \left( \sum_{j=-1}^\infty 2^{j(s-(n+1))q} \left| \sum_{\mu\in\nZ} 2^j (f^{(n+1)} , \phi_{j,\mu}) 1_{I_{j,\mu}}(\cdot) \right|^q \right)^{1/q} \right\|_{L_p} \\
    &\lesssim_s \left\| \left( \sum_{j=-1}^\infty 2^{j(s-(n+1))q} \left| \sum_{\mu\in\nZ} 2^{(j+1)(n+1)} \sum_{m\in\nZ} r_m \sfk_{j,\mu+m}(f) 1_{I_{j,\mu}}(\cdot) \right|^q \right)^{1/q} \right\|_{L_p} .
\end{align*}
Combining the powers of two and using the $u$-triangle inequality with $u = \min\{p,q,1\}$ gives
\begin{align*}
    \|f^{(n+1)}\|_{F^{s-{(n+1)}}_{p,q}} &\lesssim_s 2^{n+1} \left( \sum_{m\in\nZ} |r_m|^u \left\| \left( \sum_{j=-1}^\infty 2^{jsq} \left| \sum_{\mu\in\nZ} \sfk_{j,\mu+m}(f) 1_{I_{j,\mu}}(\cdot) \right|^q \right)^{1/q} \right\|_{L_p}^u \right)^{1/u} .
\end{align*}
Then applying Lemma \ref{lem:hardDirecMaxEstim} and finally using that $(r_m)$ decays exponentially to guarantee that the sequence $(|r_m|(|m|+1)^{1/a})_{m\in\nZ}$ is in $\ell_u$ for $a = \min\{p,q\}/2$, we see
\begin{align*}
    \|f^{(n+1)}\|_{F^{s-{(n+1)}}_{p,q}} &\lesssim_s 2^{n+1} \left( \sum_{m\in\nZ} \left(|r_m| (|m|+1)^{\frac{1}{a}}\right)^u \right)^{\frac{1}{u}}  \left\| \left( \sum_{j=-1}^\infty 2^{jsq} \left| \sum_{\mu\in\nZ} \sfk_{j,\mu}(f) 1_{I_{j,\mu}}(\cdot) \right|^q \right)^{\frac{1}{q}} \right\|_{L_p} \\
    &\lesssim_s \left\| \sfk_{j,\mu}(f) \right\|_{f^s_{p,q}} .
\end{align*}
This completes the estimate on $\|f^{(n+1)}\|_{F^{s-(n+1)}_{p,q}}$ and, therefore, the proof of the proposition.
\end{proof}

The proof of Proposition \ref{prop:hardDirecBesov} is similar and, at the same time, slightly easier.

\begin{proof}[Proof of Proposition \ref{prop:hardDirecBesov}]
Returning to \eqref{eq:BesovSumBound}, we have to estimate $\|f\|_{B^{s-(n+1)}_{p,q}}$ and $\|f^{(n+1)}\|_{B^{s-(n+1)}_{p,q}}$. We handle the former term just as in the proof of Proposition \ref{prop:hardDirecTriebel}. That is, from \cite[Section 2.3.2]{triebel1983function},
\begin{align*}
    \|f\|_{B^r_{p,q}} \lesssim_{r,\varepsilon} \|f\|_{B^{r+\varepsilon}_{p,q}} 
\end{align*}
for any $r\in\nR$ and $\varepsilon > 0$, so we can again choose $r = s-(n+1)$ and $\varepsilon = n$. Therefore, using Lemma \ref{lem:hardDirecFirstTerm} and
\begin{align*}
    \left\| (2^j (f,\psi_{j,\mu}) ) \right\|_{b^s_{p,q}} \leqslant \left\| \sfk_{j,\mu}(f) \right\|_{b^s_{p,q}}
\end{align*}
gives the desired estimate for $\|f\|_{B^{s-(n+1)}_{p,q}}$.

It remains to estimate $\|f^{(n+1)}\|_{B^{s-(n+1)}_{p,q}}$. We will write the proof as if $p,q < \infty$; the proof when $p$ or $q$ is $\infty$ can easily be written down by simply modifying each line below in the usual way. Going back to \eqref{eq:ChuiWangBesov} and using Lemma \ref{lem:hardDirecInnerBound},
\begin{align*}
    \|f^{(n+1)}\|_{B^{s-(n+1)}_{p,q}} &\lesssim_s \left( \sum_{j=-1}^{\infty} 2^{j\left( s-(n+1)-\frac{1}{p} \right)q} \left( \sum_{\mu\in\nZ} |2^j(f^{(n+1)},\phi_{j,\mu})|^p \right)^{q/p} \right)^{1/q} \\
    &\lesssim_s \left( \sum_{j=-1}^{\infty} 2^{j\left( s-(n+1)-\frac{1}{p} \right)q} \left( \sum_{\mu\in\nZ} \left| 2^{(j+1)(n+1)} \sum_{m\in\nZ} r_m \sfk_{j,\mu+m}(f) \right|^p \right)^{q/p} \right)^{1/q} .
\end{align*}
Pulling out the sum over $m$ using the $u$-triangle inequality with $u = \min\{p,q,1\}$ yields
\begin{align*}
    \|f^{(n+1)}\|_{B^{s-(n+1)}_{p,q}} &\lesssim_s 2^{n+1} \left( \sum_{m\in\nZ}  |r_m|^u \left( \sum_{j=-1}^{\infty} 2^{j\left(s-\frac{1}{p}\right)q} \left( \sum_{\mu\in\nZ} |\sfk_{j,\mu+m}(f)|^p \right)^{q/p} \right)^{u/q} \right)^{1/u} \\
    &= C_s 2^{n+1} \left( \sum_{m\in\nZ} |r_m|^u \right)^{1/u} \left( \sum_{j=-1}^{\infty} 2^{j\left(s-\frac{1}{p}\right)q} \left( \sum_{\mu\in\nZ} |\sfk_{j,\mu}(f)|^p \right)^{q/p} \right)^{1/q} \\
    &\lesssim_s \left\| \sfk_{j,\mu}(f) \right\|_{b^s_{p,q}} ,
\end{align*}
writing $C_s > 0$ to represent the implicit constant depending on $s$. Note that this time there was no trouble with the term $\sfk_{j,\mu+m}(f)$, since here we can just shift by $m$. Also, of course, the sequence $(r_m)$ is in $\ell_u$, since it decays exponentially. This completes the estimates for the final term and, therefore, the proof.
\end{proof}

Having finished the proofs for Propositions \ref{prop:hardDirecBesov} and \ref{prop:hardDirecTriebel}, we can finally give the proof for Theorems \ref{thm:mainBesov} and \ref{thm:mainTriebel}.

\begin{proof}[Proof of Theorems \ref{thm:mainBesov} and \ref{thm:mainTriebel}]
For Theorem \ref{thm:mainBesov} combine Propositions \ref{prop:easyDirecBesov} and \ref{prop:hardDirecBesov}; for Theorem \ref{thm:mainTriebel} combine Propositions \ref{prop:easyDirecTriebel} and \ref{prop:hardDirecTriebel}.
\end{proof}

\subsection{Optimality of the Lower Bound on Smoothness}\label{sec:lowerOptimal}

In this subsection, we will briefly discuss the optimality of the ranges \eqref{eq:sBesovRange} and \eqref{eq:sTriebelRange}.

\begin{prop}\label{prop:optimalityBesov}
Suppose the conditions of Theorem \ref{thm:mainBesov} and additionally that $0<p<\infty$.
\begin{enumerate}[label=(\roman*)]
    \item If $1<q<\infty$, then the lower bound for $s$ in \eqref{eq:sBesovRange} is optimal in the sense that Theorem \ref{thm:mainBesov} does not hold for $s\leqslant -n - \frac{1}{p'}$.
    \item If $0<q\leqslant 1$, then the lower bound in \eqref{eq:sBesovRange} is optimal except possibly at the endpoint, i.e. Theorem \ref{thm:mainBesov} does not hold if $s < - n - \frac{1}{p'}$.
\end{enumerate}
\end{prop}
\begin{proof}
Both of these points are simple observations about duality using \cite[Section 2.1.5]{runst2011sobolev},\footnote{We caution the reader that in this reference, when $0<p,q<1$, the dual exponents $p',q'$ are defined to be $\infty$ instead of how we have defined them in this paper.} followed by an application of Proposition \ref{prop:BLinAspq}. Consider, first, $B^s_{p,q}(\nR)$ for some $1\leqslant p<\infty$, $1<q<\infty$, and $s\leqslant -n - \frac{1}{p'}$. Then by \cite[Section 2.1.5]{runst2011sobolev}, it holds
\begin{align}
    (B^s_{p,q}(\nR))' = B^{-s}_{p',q'}(\nR) . \label{eq:BesovDual1}
\end{align}
Note, however, that $-s \geqslant n + \frac{1}{p'}$, so Proposition \ref{prop:BLinAspq} implies $\psi\notin B^{-s}_{p',q'}(\nR)$. If, on the other hand, $0<p<1$ (and still $1<q<\infty$ and $s\leqslant -n - \frac{1}{p'}$), then
\begin{align}
    (B^s_{p,q}(\nR))' = B^{-s - \frac{1}{p'}}_{\infty,q'}(\nR) , \label{eq:BesovDual2}
\end{align}
again, by \cite[Section 2.1.5]{runst2011sobolev}. By our assumption on $s$, we see that $-s - \frac{1}{p'}\geqslant n$, but by Proposition \ref{prop:BLinAspq}, we see that $\psi\in B^{-s - 1/p'}_{\infty,q'}(\nR)$ if and only if $-s - \frac{1}{p'} < n$.

The only difference in item (ii) above is that $q'$ in \eqref{eq:BesovDual1} and \eqref{eq:BesovDual2} must be replaced with $\infty$, which changes the Besov spaces that $\psi$ lives in, see Proposition \ref{prop:BLinAspq}.
\end{proof}

We see it as likely that the missing cases above are also optimal, but duality alone will not provide us with a proof of that. We are not aware of any work in this direction. Moving on to the Triebel-Lizorkin scale, much less can be said.

\begin{prop}\label{prop:optimalityTriebel}
Assume the conditions of Theorem \ref{thm:mainTriebel} and additionally $0<q<\infty$.
\begin{enumerate}[label=(\roman*)]
    \item If $1\leqslant p,q < \infty$, then Theorem \ref{thm:mainTriebel} does not hold for $s\leqslant - n - \frac{1}{p'}$.
    \item If $0<p<1$, then Theorem \ref{thm:mainTriebel} does not hold for $s < - n - \frac{1}{p'}$.
\end{enumerate}
\end{prop}
\begin{proof}
This is essentially repetitious of the proof for Proposition \ref{prop:optimalityBesov}: apply the appropriate duality results from \cite[Section 2.1.5]{runst2011sobolev} and consider in which Triebel-Lizorkin spaces $\psi$ lies using Proposition \ref{prop:BLinAspq}. We omit the details.
\end{proof}

\begin{rem}
There is a major missing piece in Proposition \ref{prop:optimalityTriebel}. Indeed, for $0<q<p<\infty$, the range \eqref{eq:sTriebelRange} becomes $-\frac{1}{q'} - n < s < n+1$. There is no immediate reason to say that the theorem could not be extended to include $-\frac{1}{p'} - n < s \leqslant -\frac{1}{q'} - n$ in this case, although it is likely, at least when $1<q<p<\infty$, that the extension is impossible. Looking at \cite[Remark 4.3]{garrigos2023haar}, the authors were able to prove that the lower bound for \eqref{eq:sTriebelRange} is optimal when $n=0$ with the help of their paper \cite[Section 6]{seeger2017haar}. Using, instead, \cite[Section 8]{srivastava2023orthogonal}, the argument in \cite[Remark 4.3]{garrigos2023haar} can be replicated in an almost identical way. Unfortunately, however, in \cite{garrigos2023haar}, the authors made use of their Proposition 4.5, which is not immediately available in our setting. Proving a version of \cite[Proposition 4.5]{garrigos2023haar} for our case is almost certainly possible and not particularly difficult, but it would require too large of an excursion to do in this paper, so we leave this question open.
\end{rem}

%% file: Sections/edgeCase.tex
\section{An Endpoint Case}\label{sec:edgeCase}

In this section, we will prove Theorem \ref{thm:edgeCase}. We will again take inspiration from \cite{garrigos2023haar}, but, also again, due to the lack of both compact support and a simple form for our wavelets, the analysis will become significantly more complicated. At any given moment in this section, all parameters (i.e. $p,q,n$) will be fixed, so dependence of constants on them will not be notated. We also assume $n\geqslant 1$. As always, $\{\psi,\Psi\}$ will represent a Battle-Lemari\'e wavelet system of order $n$. Because we are no longer using $\cB$ as our reference space in this final section, the use of Section \ref{sec:duality} in interpreting the notation $(\cdot,\cdot)$ is now unnecessary; $(\cdot,\cdot)$ can simply be interpreted as a pairing of a distribution in $\cS_{n-1}'$ with a function in $\cS_{n-1}$.


\subsection{Setup}\label{sec:edgeSetup}

Recall the material from Section \ref{sec:multiresDistribution}. With $(T_N)_{N\in\nZ}$, we will denote the multiresolution analysis of $\cS_{n-1}'$ associated to our Battle-Lemari\'e wavelet system $\{\psi,\Psi\}$. We write $\Psi_{N,\mu}(x) = \Psi(2^N x - \mu)$ for $N,\mu\in\nZ$ and define the projection onto $T_N$: for $N\in\nZ$,
\begin{align}
    E_N f = \sum_{\mu\in\nZ} 2^N (f,\Psi_{N,\mu}) \Psi_{N,\mu} \label{eq:projection}
\end{align}
for $f\in \cS_{n-1}'(\nR)$. As we mentioned in Section \ref{sec:multiresDistribution}, we know $E_N f\to f$ in $\cS_{n-1}'$ for $f\in\cS_{n-1}'$.


The other definition that we will need is an antiderivative for our Battle-Lemari\'e wavelet, $\psi$. We will define $\rho$ so that it satisfies $\rho^{(n+1)}(2x) = \psi(x)$. Of course, this does not uniquely determine $\rho$, since any polynomial of order $n$ or less will be destroyed by the $n+1$ derivatives. In Lemma \ref{lem:rhoExpDecay} below, we will show that we can choose $\rho$ to have exponential decay.

\begin{lem}\label{lem:rhoExpDecay}
Let $n\in\nN_0$ and suppose $\psi$ is a Battle-Lemari\'e wavelet of order $n$. The function $\rho$ defined through $\rho^{(n+1)}(2x) = \psi(x)$ can be chosen to decay exponentially. Additionally, when $n=0$, $\rho$ can be chosen to have compact support.
\end{lem}

In order to prove this, we give a nod to G. Battle in \cite{battle1987block}, who proved that the original Battle-Lemari\'e scaling function decays exponentially using methods from complex analysis, specifically Paley-Wiener theorems. Before beginning the proof of the above lemma, we collect the results from complex analysis that we will need for the reader's convenience. We will follow \cite[Chapter 4]{stein2010complex}. Define $\mathfrak{F}_a$ for $a>0$ as those complex-valued functions, $f$, which are holomorphic in $S_a = \{z\in\nC : |\text{Im}(z)| < a\}$ and which satisfy
\begin{align}
    |f(x + iy)| \leqslant \frac{C}{1 + x^2} \label{eq:FaBound}
\end{align}
for all $x\in\nR$, $|y|<a$, and some constant $C > 0$.

\begin{prop}[Theorems 2.1 and 3.1 in Chapter 4 of \cite{stein2010complex}]\label{prop:PaleyWiener}
Let $a > 0$. If $f\in\mathfrak{F}_a$, then $\widehat{f}$ decays exponentially. Conversely, suppose $f$ is such that $|f(x)|\lesssim 1/(1+x^2)$ for all $x\in\nR$ and $|\widehat{f}(\xi)| \leqslant B e^{-b|\xi|}$ for $\xi\in\nR$ and some constants $b,B > 0$. Then, for any $0<c<b$, we have that $f$ is the restriction to $\nR$ of a function, which is holomorphic in $S_c$.
\end{prop}

\begin{rem}\label{rem:extUniformBound}
Checking the proof of \cite[Theorem 3.1 of Chapter 4]{stein2010complex} (using the notation and assumptions of Proposition \ref{prop:PaleyWiener}), it should come as no surprise that the extension of $f$ to $S_c$ is given by
\begin{align*}
    \int_{-\infty}^\infty \widehat{f}(\xi) e^{i\xi z} \, d\xi , \quad z\in S_c .
\end{align*}
Thus, by the exponential decay of $\widehat{f}$, this extension is uniformly bounded for $z\in S_c$ by
\begin{align*}
    B \int_{-\infty}^\infty e^{-b|\xi|} e^{c|\xi|} \, d\xi ,
\end{align*}
which is why we had to assume $0<c<b$.
\end{rem}



With this reminder in hand, we are ready to prove the lemma.

\begin{proof}[Proof of Lemma \ref{lem:rhoExpDecay}]
We begin by distinguishing the case $n=0$, which can be proven essentially by inspection. The classical Battle-Lemari\'e wavelet of order zero is the Haar function, see \eqref{eq:HaarSystem}. We have already seen in \eqref{eq:Haar2scale} that this lemma holds in that case. Further, looking back at Remark \ref{rem:Haar}, any other Battle-Lemari\'e wavelet of order zero, call it $\psi$, satisfies $\psi(x) = h(x-k)$ for some $k\in\nZ$, so the lemma holds in that case as well.

For the remainder of the proof, $n\in\nN$ is fixed and so is the Battle-Lemari\'e wavelet, $\psi$, of order $n$. Considering Proposition \ref{prop:PaleyWiener}, our goal shall be to prove that we can extend $\widehat{\rho}$ to be an element of $\mathfrak{F}_a$ for some $a>0$, preliminarily call this extension $P$ (we will be more precise momentarily). The proposition will, then, imply that $\widehat{P}$ decays exponentially, implying that its restriction, $\rho$, does as well ($P$ restricted to $\nR$ is $\widehat{\rho}$, so $\widehat{P}$ restricted to $\nR$ is $\widehat{\widehat{\rho}} = \rho(-\cdot)$).

Let us make some observations. First, we know that the Fourier transform of $\widehat{\psi}$, i.e. $\widehat{\widehat{\psi}}$, decays exponentially; hence, by the converse part in Proposition \ref{prop:PaleyWiener}, we conclude that $\widehat{\psi}$ is the restriction to $\nR$ of a function, $Q(z)$, holomorphic in the horizontal strip $S_{2a}$ for some $a > 0$.

Second, based on the requirement $\rho^{(n+1)}(2x) = \psi(x)$,
\begin{align*}
    \widehat{\rho}(\xi) = \frac{2}{(i\xi)^{n+1}} \widehat{\psi}(2\xi) , \quad \xi\in\nR\smallsetminus\{0\} .
\end{align*}
We already know that $\widehat{\psi}(2\xi)$ can be extended holomorphically to $S_a$ and clearly $\frac{2}{(i\xi)^{n+1}}$ can as well except at $\xi = 0$. Write $\widetilde{P}$ for this extension of $\widehat{\rho}$ to $S_a\smallsetminus \{0\}$.

In the terminology of complex analysis, this singularity at the origin is called removable: note that
\begin{align}
    \lim_{z\to 0} z \frac{1}{z^{n+1}} Q(2z) = 0 . \label{eq:riemannLimit}
\end{align}
Indeed, if we expand $Q$ as a power series around zero, using that $Q^{(\kappa)}(0) = 0$ for $0\leqslant \kappa\leqslant n$ because of the vanishing moments of $\psi$, we can swiftly prove that the above limit is null. Hence, by Riemann's famous theorem on removable singularities, see \cite[Theorem 7 in Chapter 4]{martin1966complex},\footnote{Roughly, the theorem says that if the limit in \eqref{eq:riemannLimit} is zero, then the singularity is removable, i.e. there exists a unique extension.} $\widetilde{P}$ can be extended uniquely to a holomorphic function on the whole of $S_a$; we will call this extension $P$.

Returning to Proposition \ref{prop:PaleyWiener}, if we can establish \eqref{eq:FaBound} for $P$, then we are done. In fact, we have already established that $P(z)$ is bounded near zero (it has a removable singularity there). Thus, it suffices to show \eqref{eq:FaBound} away from the origin. Let $U$ be some small neighborhood of the origin; for concreteness, say $U = \{z\in \nC : |z| < 1\}$. For $z\in S_a \smallsetminus U$,
\begin{align*}
    |P(z)| = \frac{2}{|z|^{n+1}} |Q(2z)| \leqslant \frac{4}{1 + |z|^{n+1}} |Q(2z)| \leqslant \frac{4}{1+|\text{Re}(z)|^{n+1}} |Q(2z)| .
\end{align*}
The second inequality follows from the elementary estimate $|\text{Re}(z)| \leqslant |z|$ for all $z\in\nC$. Now, by Remark \ref{rem:extUniformBound}, the extension $Q$ is uniformly bounded, say by $Q_0 > 0$. Hence,
\begin{align*}
    |P(x + iy)| \leqslant \frac{4 Q_0}{1 + |x|^{n+1}}
\end{align*}
for $x+iy = z$, where $z\in S_a\smallsetminus U$. Since $n\geqslant 1$, we have established \eqref{eq:FaBound} for $P$, so we are done.
\end{proof}


\subsection{The Weaker Direction}\label{sec:edgeCaseWeaker}

With the setup out of the way, we proceed with the proof of Theorem \ref{thm:edgeCase}. First, we restate and prove Proposition \ref{prop:p1inftyEmbed}.

\begin{prop}[Restatement of Proposition \ref{prop:p1inftyEmbed}]\label{prop:forwardEdge}
Suppose $n\geqslant 1$, $f\in \cS_{n-1}'(\nR)$, and $\{\psi,\Psi\}$ is a Battle-Lemari\'e system of order $n$. Let $1\leqslant p\leqslant\infty$. There exists a constant $C = C(p,n,\psi,\Psi) > 0$ such that
\begin{align*}
    \sup_{j\geqslant -1} 2^{j\left(n+1 - \frac{1}{p}\right)} \left( \sum_{\mu\in\nZ} |\sfk_{j,\mu}(f)|^p \right)^{1/p} \leqslant C\|f\|_{W^{n+1}_p(\nR)} ,
\end{align*}
making the usual modification if $p=\infty$.
\end{prop}
\begin{proof}
Begin by noting that since $\psi,\Psi\in \cS_{n-1}$, the shifted spline coefficients, $\sfk_{j,\mu}(f)$ make sense for all $f\in \cS_{n-1}'(\nR)$ (and $j\geqslant -1, \mu\in\nZ$). For the remainder of the proof, we may assume that $f\in W^{n+1}_p(\nR)$; otherwise the desired inequality is vacuous.

We define $\rho$ as in Lemma \ref{lem:rhoExpDecay} and choose it so that it decays exponentially. Note that, by definition, it holds for $f\in W^{n+1}_p(\nR)$ that
\begin{align}
    (f,\widetilde{\psi}_{j-1,\nu}) = (-1)^{n+1} 2^{-j(n+1)} (f^{(n+1)},\rho_{j,\nu})  \label{eq:rhoLemma2.1}
\end{align}
for $j\in\nN$ and $\nu\in\nZ$. In fact, we have defined $\rho$ with the intention of producing \eqref{eq:rhoLemma2.1} (it is made to resemble \cite[Lemma 2.1]{garrigos2023haar}).

We fix $j\geqslant 0$ and $\nu\in\nZ$; due to \eqref{eq:rhoLemma2.1},
\begin{align}
    |2^{(j+1)(n+1)} (f,\widetilde{\psi}_{j,\nu})| &\leqslant \int |f^{(n+1)} \rho_{j+1,\nu}| \, dx = \sum_{\mu\in\nZ} \int_{I_{j+1,\mu}} |f^{(n+1)}\rho_{j+1,\nu}| \, dx \notag \\
    &\leqslant \sum_{\mu\in\nZ} \|f^{(n+1)}\|_{L_p(I_{j+1,\mu})} \|\rho_{j+1,\nu}\|_{L_{p'}(I_{j+1,\mu})} \label{eq:step4.1} .
\end{align}
Recall that by Lemma \ref{lem:rhoExpDecay}, $\rho$ decays exponentially, i.e. $|\rho(x)|\lesssim e^{-\delta|x|}$ for some $\delta > 0$. If $1\leqslant p' < \infty$, the norm of $\rho$ can be rewritten as
\begin{align*}
    \|\rho_{j+1,\nu}\|_{L_{p'}(I_{j+1,\mu})}^{p'} = 2^{-(j+1)} \int_{\mu-\nu}^{\mu-\nu+1} |\rho(x)|^{p'} \, dx \lesssim 2^{-(j+1)} e^{-\delta p' |\mu-\nu|} .
\end{align*}
When $p' = \infty$, we also get $\|\rho_{j+1,\nu}\|_{L_\infty(I_{j+1,\mu})} \lesssim e^{-\delta |\mu-\nu|}$. Inserting this back into \eqref{eq:step4.1},
\begin{align*}
    |2^{(j+1)(n+1)} (f,\widetilde{\psi}_{j,\nu})| \lesssim \sum_{\mu\in\nZ} 2^{-(j+1)/p'} e^{-\delta|\mu-\nu|} \|f^{(n+1)}\|_{L_p(I_{j+1,\mu})} .
\end{align*}
Hence, assuming $p<\infty$,
\begin{align*}
    &\sup_{j\geqslant 0} 2^{j\left( n+1 - \frac{1}{p} \right)} \left( \sum_{\nu\in\nZ} |2^j(f,\widetilde{\psi}_{j,\nu})|^p \right)^{1/p} \\
    &\hspace{+5mm}\lesssim \sup_{j\geqslant 0} 2^{j\left( n+1 - \frac{1}{p} \right)} 2^{-j(n+1)} 2^{-(n+1)} 2^{-(j+1)\left( 1 - \frac{1}{p} \right)} 2^j  \left( \sum_{\nu\in\nZ} \left| \sum_{\mu\in\nZ} e^{-\delta|\mu-\nu|} \|f^{(n+1)}\|_{L_p(I_{j+1,\mu})} \right|^p \right)^{1/p} .
\end{align*}
The powers of $2$ all simplify to $2^{-(n+1)} 2^{-1/p'}$, which we bound with $1$. Shifting $\mu-\nu\mapsto \mu$ on the inside sum, this becomes
\begin{align*}
    \sup_{j\geqslant 0} \left( \sum_{\nu\in\nZ} \left| \sum_{\mu\in\nZ} e^{-\delta|\mu|} \|f^{(n+1)}\|_{L_p(I_{j+1,\mu+\nu})} \right|^p \right)^{1/p} &\leqslant \sup_{j\geqslant 0} \sum_{\mu\in\nZ} \left( \sum_{\nu\in\nZ} \left| e^{-\delta |\mu|} \|f^{(n+1)}\|_{L_p(I_{j+1,\mu+\nu})} \right|^p \right)^{1/p} \\
    &= \sup_{j\geqslant 0} \left( \sum_{\mu\in\nZ} e^{-\delta |\mu|} \right) \left( \sum_{\nu\in\nZ} \|f^{(n+1)}\|_{L_p(I_{j+1,\nu})}^p \right)^{1/p} \\
    &= C_\delta \|f^{(n+1)}\|_{L_p(\nR)} .
\end{align*}
The inequality above is just the triangle inequality for $\ell_p$; we have also translated $\mu+\nu\mapsto \nu$ in the sum over $\nu$ in the penultimate step. The above steps can be repeated when $p=\infty$, just making the obvious modifications; the result is the same.

In the case $j = -1$, we can write for $\nu\in\nZ$,
\begin{align*}
    |(f,\widetilde{\psi}_{-1,\nu})| &= \left| \int f(x) \Psi(x-\nu) \, dx \right| \leqslant \sum_{\mu\in\nZ} \int_{I_{-1,\mu}} |f(x)\Psi(x-\nu)| \, dx \\
    &\leqslant \sum_{\mu\in\nZ} \|f\|_{L_p(I_{-1,\mu})} \|\Psi(\cdot - \nu)\|_{L_{p'}(I_{-1,\mu})} .
\end{align*}
Trivially, $\|\Psi(\cdot - \nu)\|_{L_{p'}(I_{-1,\mu})} \lesssim e^{-\gamma |\mu-\nu|}$. Therefore, if $p<\infty$, taking similar steps as before,
\begin{align*}
    2^{-\left( n+1-\frac{1}{p} \right)} \left( \sum_{\nu\in\nZ} |(f,\widetilde{\psi}_{-1,\nu})|^p \right)^{1/p} &\lesssim \left( \sum_{\nu\in\nZ}  \left| \sum_{\mu\in\nZ} e^{-\gamma |\mu-\nu|} \|f\|_{L_p(I_{-1,\mu})} \right|^p \right)^{1/p} \\
    &\leqslant \sum_{\mu\in\nZ} \left( e^{-\gamma|\mu|} \right) \left( \sum_{\nu\in\nZ} \|f\|_{L_p(I_{-1,\nu})}^p \right)^{1/p} \\
    &= C_\gamma \|f\|_{L_p(\nR)} .
\end{align*}
If $p=\infty$, we can once again just repeat the above steps, which completes the proof.
\end{proof}


\subsection{The Stronger Direction}\label{sec:edgeStronger}

This brings us to the final loose end.

\begin{prop}\label{prop:backwardEdge}
Suppose $n\geqslant 1$, $f\in \cS_{n-1}'(\nR)$, and $\{\psi,\Psi\}$ is a Battle-Lemari\'e system of order $n$. Let $1<p\leqslant \infty$. There exists a constant $C = C(p,n,\psi,\Psi) > 0$ such that
\begin{align*}
    \|f\|_{W^{n+1}_p(\nR)}\leqslant C \sup_{j\geqslant -1} 2^{j\left(n+1 - \frac{1}{p}\right)} \left( \sum_{\mu\in\nZ} |\sfk_{j,\mu}(f)|^p \right)^{1/p} ,
\end{align*}
making the usual modifications if $p=\infty$.
\end{prop}
\begin{proof}
We start by mentioning that, as we have explained at the beginning of the proof for Proposition \ref{prop:forwardEdge}, the shifted spline coefficients make sense for $f\in\cS_{n-1}'$, because $\psi,\Psi\in\cS_{n-1}$.

Our first goal is to show
\begin{align*}
    \|f\|_{L_p} \lesssim \sup_{j\geqslant -1} 2^{j\left( n + 1 - \frac{1}{p} \right)} \left( \sum_{\mu\in\nZ} |2^{j_+} (f,\widetilde{\psi}_{j,\mu})|^p \right)^{1/p} =: A_0 
\end{align*}
for $1 \leqslant p \leqslant \infty$. We will start by establishing
\begin{align*}
    \sup_{N\geqslant 0} \|E_N f\|_{L_p} \lesssim A_0 .
\end{align*}
Recall the work in Section \ref{sec:multiresDistribution}. Applying $T_N = T_0 \dotplus U_0 \dotplus \cdots \dotplus U_{N-1}$ for $N\in\nN$, we can decompose $E_N f$ as
\begin{align}
    E_N f = E_0 f + \sum_{0\leqslant j < N} \sum_{\mu\in\nZ} 2^j (f,\psi_{j,\mu}) \psi_{j,\mu} \label{eq:ENdecomp}
\end{align}
for all $f\in\cS_{n-1}'(\nR)$.

We look at the sum over $\mu$, first when $j\geqslant 0$. Observe,
\begin{align*}
    \left\| \sum_{\mu\in\nZ} 2^j (f,\psi_{j,\mu}) \psi_{j,\mu} \right\|_{L_p} \lesssim 2^{-j/p} \left( \sum_{\mu\in\nZ} |2^j (f,\psi_{j,\mu})|^p \right)^{1/p} \leqslant 2^{-j(n+1)} A_0 .
\end{align*}
The first inequality is not immediately obvious for $1<p<\infty$, but is a standard fact, see \cite[Lemma 8.2]{wojtaszczyk1997mathematical}.\footnote{We could also justify this step with an unnecessarily big hammer: look at $L_p$ as $F^0_{p,2}$ for $1<p<\infty$ and use Theorem \ref{thm:SplineCharac}. It is, then, obvious, because the dyadic intervals, $I_{j,\mu}$, are disjoint for fixed $j$.} We can repeat the steps above with $j=-1$; this yields $\|E_0 f\|_{L_p}\lesssim A_0$. Furthermore, if $N_2 < N_1$ are nonnegative integers, then using \eqref{eq:ENdecomp},
\begin{align*}
    \|E_{N_1} f - E_{N_2} f\|_{L_p} &= \left\| \sum_{N_2\leqslant j < N_1} \sum_{\mu\in\nZ} 2^j (f,\psi_{j,\mu})\psi_{j,\mu} \right\|_{L_p} \lesssim \sum_{N_2\leqslant j < N_1} 2^{-j(n+1)} A_0 \\
    &\leqslant 2\cdot 2^{-N_2} A_0 .
\end{align*}
This proves that $E_N f$ converges in $L_p$, although it does not say what the limit function is. We mentioned in Section \ref{sec:multiresDistribution} that, for $f\in\cS_{n-1}'$, we have $E_N f\to f$ in $\cS_{n-1}'$ as $N\to\infty$, which is weaker than convergence in $L_p$, so $E_N f\to f$ in $L_p$. It also shows
\begin{align*}
    \sup_{N\geqslant 0} \|E_N f\|_{L_p} \leqslant \|E_0 f\|_{L_p} + \sup_{N\geqslant 0} \|E_N f - E_0 f\|_{L_p} \lesssim A_0 ,
\end{align*}
completing the proof of the claim. With this observation, for all $N\in\nN_0$
\begin{align*}
    \|f\|_{L_p} \leqslant \|f - E_N f\|_{L_p} + \sup_{N\geqslant 0} \|E_N f\|_{L_p} ,
\end{align*}
which verifies that $\|f\|_{L_p} \lesssim A_0$.

We now finish this proof by showing
\begin{align*}
    \|f^{(n+1)}\|_{L_p} \lesssim \sup_{j\geqslant -1} 2^{j\left( n+1 - \frac{1}{p} \right)} \left( \sum_{\mu\in\nZ} |2^{j_+} (f,\widetilde{\psi}_{j,\mu})|^p \right)^{1/p} = A_0
\end{align*}
for $1<p\leqslant \infty$. As in Section \ref{sec:strong}, we will take $\{\phi,\Phi\}$ to be one of the Battle-Lemari\'e wavelet systems of order $2n+1$ from \cite{ushakova2018localisation}. We will write $\Phi_{N,\mu}(x) = \Phi(2^N x - \mu)$ for $N,\mu\in\nZ$. Recall from Section \ref{subsubsec:remarkableConnection} that there exists a sequence $(c_k)\subset\nR$ decaying exponentially with
\begin{align*}
    \Phi(x) = \sum_{k\in\nZ} c_k B_{2n+1}(x-k) .
\end{align*}
Furthermore, using similar steps to those in the proof of Lemma \ref{lem:hardDirecInnerBound}, we note that for $N\in\nN$ and $\mu\in\nZ$,
\begin{align}
    \Phi^{(n+1)}(2^{N+1}x-\mu) = \sum_{\ell\in\nZ} t_\ell \widetilde{\psi}_{N,\ell+\mu} , \label{eq:PhiDerivative}
\end{align}
where $t_\ell$ is given by
\begin{align*}
    t_\ell = \sum_{k\in\nZ} c_k q_{\ell-k}
\end{align*}
for each $\ell\in\nZ$. Recall the sequence $(q_\ell)$ comes from Lemma \ref{lem:hardPartCrux} and decays exponentially. In particular, $(t_\ell)\subset\nR$ decays exponentially too by \cite[Lemma 3.16]{wojtaszczyk1997mathematical}. We also define the projection onto $V_N$, the scale space for $N\in\nZ$ corresponding to $\Phi$:
\begin{align*}
    P_N g(x) = \sum_{\mu\in\nZ} 2^N (g,\Phi_{N,\mu}) \Phi_{N,\mu}(x)
\end{align*}
for $g\in L_2(\nR)$. Fix $N\in\nN$ and define
\begin{align*}
    g_N(x) = \sum_{\mu\in\nZ} 2^{N+1} (f^{(n+1)}, \Phi_{N+1,\mu}) \Phi_{N+1,\mu}(x) .
\end{align*}
Of course, $(f^{(n+1)}, \Phi_{N+1,\mu}) = (-1)^{N+1} 2^{(N+1)(n+1)} (f,\Phi_{N+1,\mu}^{(n+1)})$. Thus, using \eqref{eq:PhiDerivative},
\begin{align*}
    \left| 2^{N+1} (f^{(n+1)},\Phi_{N+1,\mu}) \right| \leqslant 2^{N+1} 2^{(N+1)(n+1)} \sum_{\ell\in\nZ} |t_\ell| |(f,\widetilde{\psi}_{N,\ell+\mu})| .
\end{align*}
Below we will write the proof for $p<\infty$; the analogous proof works for $p=\infty$. Using, again, \cite[Lemma 8.2]{wojtaszczyk1997mathematical} we see
\begin{align*}
    \|g_N\|_{L_p} &= \left\| \sum_{\mu\in\nZ} 2^{N+1} (f^{(n+1)}, \Phi_{N+1,\mu}) \Phi_{N+1,\mu} \right\|_{L_p} \\
    &\lesssim 2^{-(N+1)/p} \left( \sum_{\mu\in\nZ} |2^{N+1} (f^{(n+1)}, \Phi_{N+1,\mu})|^p \right)^{1/p} .
\end{align*}
Hence,
\begin{align*}
    \|g_N\|_{L_p} &\lesssim 2^{-(N+1)/p} \left( \sum_{\mu\in\nZ} \left| 2^{N+1} 2^{(N+1)(n+1)} \sum_{\ell\in\nZ} |t_\ell| |(f,\widetilde{\psi}_{N,\ell+\mu})| \right|^p \right)^{1/p} \\
    &\leqslant 2\cdot 2^{n+1}2^{-1/p} 2^{N\left( n+1-\frac{1}{p} \right)} \sum_{\ell\in\nZ}  \left( |t_\ell| \left[ \sum_{\mu\in\nZ} |2^N (f,\widetilde{\psi}_{N,\ell+\mu})|^p \right]^{1/p} \right) \\
    &\lesssim 2^{N\left( n+1-\frac{1}{p}\right)} \left( \sum_{\mu\in\nZ} |2^N (f,\widetilde{\psi}_{N,\mu})|^p \right)^{1/p} \\
    &\leqslant A_0 .
\end{align*}
In the second line above, we just used the $\ell_p$ triangle inequality.

Using the (sequential) Banach-Alaoglu theorem, we have, therefore, proven, that there exists $g\in L_p(\nR)$ that is the weak-$*$ limit of a subsequence of $(g_N)$. Noting that for $N+1\geqslant j$, we have $P_{N+1}(\phi_{j,\nu}) = \phi_{j,\nu}$ for $\nu\in\nZ$, it follows that
\begin{align*}
    (g_N , \phi_{j,\nu}) = \left( \sum_{\mu\in\nZ} 2^{N+1} (f^{(n+1)} , \Phi_{N+1,\mu}) \Phi_{N+1,\mu} , \phi_{j,\nu} \right) = (f^{(n+1)} , \phi_{j,\nu}) 
\end{align*}
for all $-1\leqslant j\leqslant N+1$ and $\nu\in\nZ$. Taking the limit as $N\to\infty$, we see $(g,\phi_{j,\nu}) = (f^{(n+1)}, \phi_{j,\nu})$ for $j\geqslant -1$ and $\nu\in\nZ$. Since $\{\phi,\Phi\}$ gives an orthogonal basis of $L_2$, any Schwartz function can be decomposed with them. This, then, implies that $g = f^{(n+1)}$ as tempered distributions and therefore as $L_p$ functions (since $g\in L_p$). Hence, 
\begin{align*}
    \|f^{(n+1)}\|_{L_p} = \|g\|_{L_p}\lesssim A_0 .
\end{align*}
This completes the proof. 
\end{proof}

We can now prove Theorem \ref{thm:edgeCase}.

\begin{proof}[Proof of Theorem \ref{thm:edgeCase}]
Combine Propositions \ref{prop:forwardEdge} and \ref{prop:backwardEdge}. 
\end{proof}

\begin{rem}\label{rem:BV}
When $p=1$, the entire proof of this direction works until we apply the Banach-Alaoglu theorem, which can still be applied, it is just that we would only find that $f^{(n+1)}$ has total variation bounded by $A_0$. The argument can be repeated for any derivative of $f$ up to order $n+1$ to see that each derivative has total variation bounded by $A_0$. This may be of some interest, but we do not pursue this line of reasoning here.
\end{rem}